\documentclass[]{article}

\usepackage[utf8]{inputenc}
\usepackage[T1]{fontenc}

\usepackage[hmarginratio=1:1,top=32mm,columnsep=20pt]{geometry} % Document margins
\usepackage[hang, small,labelfont=bf,up,textfont=it,up]{caption} % Custom captions under/above floats in tables or figures
\usepackage{hyperref}
\usepackage{siunitx}

\usepackage{amsmath}
\usepackage{mathtools}
\usepackage{cases}
\newcommand{\rom}{\ensuremath{\text{ROM}}}
\newcommand{\podec}{\ensuremath{\text{POD}}}
\newcommand{\finca}{\ensuremath{\text{Fincantieri S.p.A.}}}

\newcommand{\cfd}{\ensuremath{\text{CFD}}}
\newcommand{\fom}{\ensuremath{\text{FOM}}}
\newcommand{\nse}{\ensuremath{\text{NSE}}}
\newcommand{\mrf}{\ensuremath{\text{MRF}}}
\newcommand{\rbf}{\ensuremath{\text{RBF}}}
\newcommand{\gpr}{\ensuremath{\text{GPR}}}
\newcommand{\knr}{\ensuremath{\text{KNR}}}
\newcommand{\nsnaps}{\ensuremath{\text{M}}}
\newcommand{\nparams}{\ensuremath{\text{p}}}

\usepackage{booktabs}
\newcommand{\RA}[1]{{\color{black}#1}}
\newcommand{\RB}[1]{{\color{black}#1}}

\usepackage{subfig}
\usepackage{multirow}
\newcolumntype{M}[1]{>{\centering\arraybackslash}m{#1}}

\usepackage{comment}

\usepackage{amssymb, latexsym, amsmath}
\usepackage{xcolor}
\usepackage{colortbl}

\usepackage{titlesec}
\titleformat{\section}[block]{\large\scshape\centering}{\thesection.}{1em}{} % Change the look of the section titles
\titleformat{\subsection}[block]{\large}{\thesubsection.}{1em}{} % Change the look of the section titles

\usepackage{abstract} % Allows abstract customization
\title{\textbf{A shape optimization pipeline for marine propellers by means of reduced order modeling techniques}}

\date{ }
\author{Anna Ivagnes  \\ \small SISSA, International School for Advanced Studies, \\ \small Mathematics Area, mathLab, Trieste, Italy. \\ \small  \href{mailto:aivagnes@sissa.it}{aivagnes@sissa.it} \normalsize \and Nicola Demo \\ \small SISSA, International School for Advanced Studies, \\ \small Mathematics Area, mathLab, Trieste, Italy. \\ \small  \href{mailto:ndemo@sissa.it}{ndemo@sissa.it} \normalsize  \and Gianluigi Rozza \\ \small  SISSA, International School for Advanced Studies, \\ \small Mathematics Area, mathLab, Trieste, Italy. \\ \small \href{mailto:grozza@sissa.it}{grozza@sissa.it} }
\begin{document}
\maketitle
\begin{abstract}
\noindent  In this paper, we propose a shape optimization pipeline for propeller blades, applied to naval applications. The geometrical features of a blade are exploited to parametrize it, allowing to obtain deformed blades by perturbating their parameters. The optimization is performed using a genetic algorithm that exploits the computational speed-up of reduced order models to maximize the efficiency of a given propeller. A standard offline-online procedure is exploited to construct the reduced-order model. In an expensive offline phase, the full order model, which reproduces an open water test, is set up in the open-source software OpenFOAM and the same full order setting is used to run the CFD simulations for all the deformed propellers. The collected high-fidelity snapshots and the deformed parameters are used in the online stage to build the non-intrusive reduced-order model. This paper provides a proof of concept of the pipeline proposed, where the optimized propeller improves the efficiency of the original propeller. 
\end{abstract}

\section{Introduction}
\label{sec:intro}
The problem of optimizing the efficiency of marine propellers is a topic of considerable importance in naval engineering applications.
The correct understanding of how the geometry of the blades would impact the propagation of vibrations and noise is propaedeutic to the design of new propellers and for the improvement of the propulsion performances.
Therefore, a preliminary but necessary step in the optimization of propellers' shape is the geometrical parametrization of a single blade, where both the global geometrical features, i.e. pitch, rake, skew, chord lengths, and the sections' properties, i.e. camber, thickness, are recognized as parameters~\cite{carlton2018marine, kerwin1986marine}.
Modifications in the parameters lead to different blade shapes, allowing the exploration of a large number of different shapes.
%In particular, in our analysis, we will modify the chord lengths, the pitch, the camber and the thickness.
In the project here presented, we started from a blade provided by \finca{}, but the starting blade can be in general defined starting from its parameters.\\

The step following the parametrization is the setup of the Computational Fluid Dynamic (\cfd{}) simulation in order to compute the propeller performances.
In particular, different simulation settings have been investigated in the state of the art: \emph{open water} tests~\cite{baltazar2012open,boswell1971design,bhattacharyya2016scale}, usually studied at a reduced scale, and the inverse problem concerning the wake generated by the propeller, analyzed in \RB{previous works, such as}~\cite{muscari2013modeling,zhang2019numerical}.
In particular, this paper focuses on the optimization of propeller efficiency in open-water tests.
%and the setting of experimental tests is reproduced by using the open-source finite-volume based software \emph{OpenFOAM}.
From a computational point of view, the main challenge is to find a compromise between the computational cost of high-fidelity simulations, i.e. the mesh refinement, the computational resources, and the precision in the reconstruction of the most significant physical fields in experimental tests, i.e. the torque and the thrust coefficients. Indeed, an accurate reconstruction of the experimental tests would require a significantly large number of degrees of freedom and the computational effort needed for the optimization process would become unfeasible.\\

%In particular, in the project here presented, the computational mesh is composed of $\simeq \num{6000000}$ cells, which provides an error lower than $\simeq 3 \%$ with respect to the quantities of interest.
To mitigate this issue, this contribution introduces an optimization pipeline exploiting the potentialities of \emph{data-driven} or \emph{non-intrusive} Reduced Order Models (\rom{}s) \cite{salmoiraghi2018free,quarteroni2015reduced,rozza2015book,morhandbook2020, rozza2022advanced, degruyter1, degruyter2, degruyter3}. \rom{}s constitutes a consolidated method for real-time approximation of the numerical solutions in different problem configurations, i.e. finite volume, finite element, and so on, and in naval applications. We can find evidence of the \rom{} potentiality in naval problems in previous works, such that~\cite{demo2021hull, demoortaligustinrozzalavini2020bumi, DemoTezzeleMolaRozza2021JMSE, tezzele2018dimension, tezzele2020enhancing, tezzele2022integrated, tezzele2023multifidelity, d2018combined, ccelik2021reduced, serani2019stochastic}.

The standard approach when dealing with \rom{}s is the \emph{offline-online} procedure. The offline stage consists of the computation of a large number of high-fidelity simulations, each one exploiting a different deformed blade; the full-order snapshots of the main flow fields, i.e. pressure and wall shear stress are then collected. In the online stage, the reduced order model built from snapshots and parameters is exploited to predict the propeller efficiency in the optimization process.

This contribution proposes two different \rom{}s approaches: a standard one, where the snapshots are evaluated on all the mesh points of the blades, and a \emph{fast} \rom{}, where the fields are evaluated only on a reduced number of blades points, i.e. the Gauss quadrature nodes. To the knowledge of the authors, this approach is here experimented for the first time in a \rom{} fashion and allows for a significant \RB{reduction in the computational time} without losing the \rom{} accuracy.

\begin{figure*}
    \centering
    \includegraphics[width=0.8\textwidth]{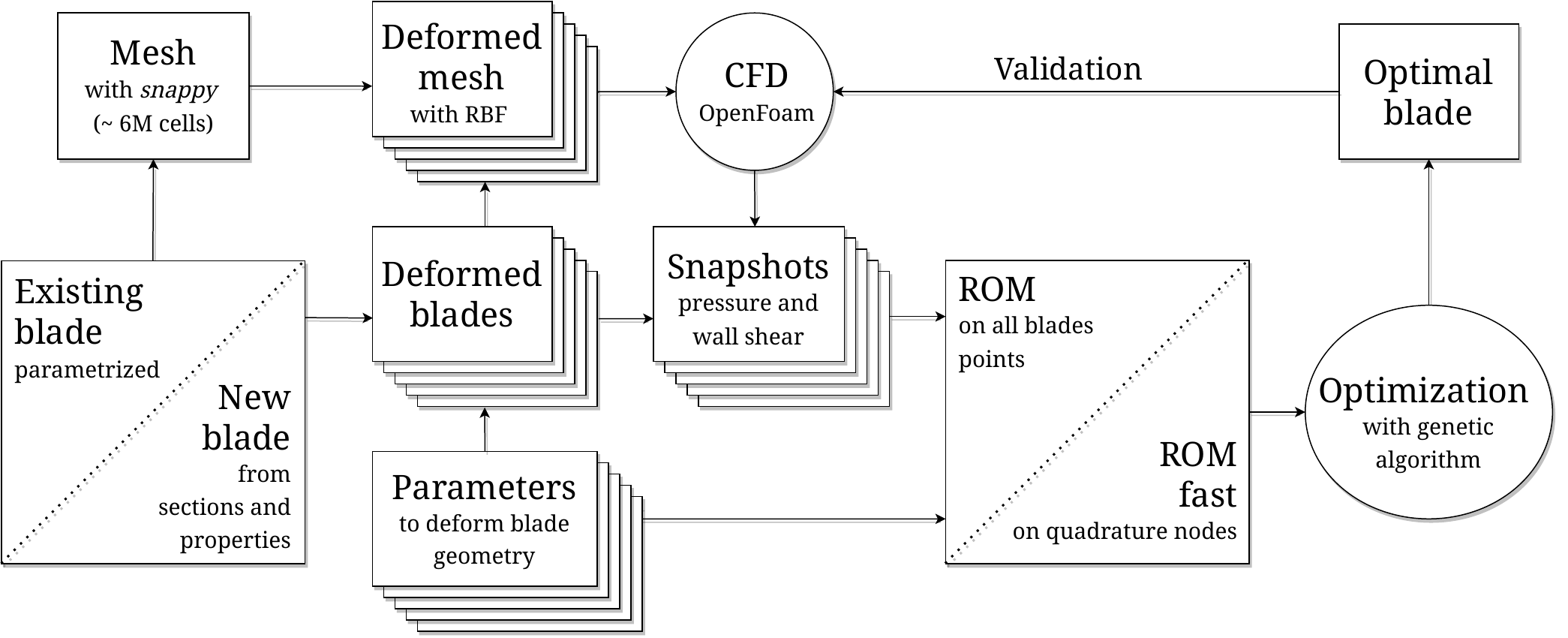}
    \caption{Pipeline of the shape optimization.}
    \label{fig:pipeline}
\end{figure*}

We now present the structure of the paper, which follows the pipeline represented in Figure \ref{fig:pipeline}. Section \ref{sec:param-def} will analyze the geometry of a propeller and of a single blade (Subsection \ref{subsec:param}), the deformation of the original propeller (Subsection \ref{subsec:blade-def}), and of the whole mesh (Subsection \ref{subsec:mesh-def}). The setting of the Full Order Model (\fom{}), which is simulated in the open-source software \emph{OpenFOAM}, is described in detail in Subsection \ref{sec:fom}. Then, the basic theory of Reduced Order Models (\rom{}) is presented in Section \ref{sec:roms}, where two different \rom{}s approaches are described (Subsections \ref{rom_1} and \ref{rom_2}). Section \ref{sec:opt} describes how the two \rom{}s techniques are exploited either in a genetic (Subsection \ref{subsec:genetic}) or in a gradient-based method (Subsection \ref{subsec:gradient}). The results obtained in all the optimization processes are there compared.

\section{Propeller geometry and mesh deformation}
\label{sec:param-def}
This Section is dedicated to the analysis of the geometry of a single blade, followed by the exploitation of its main parameters to obtain different deformed shapes (Subsections \ref{subsec:param} and \ref{subsec:blade-def}).

The deformation of the blade would result in the deformation of the whole \emph{OpenFOAM} mesh in our offline simulations. This specific issue is addressed in Subsection \ref{subsec:mesh-def}.

\subsection{Blade parameters}
\label{subsec:param}

The first important step in the development of efficient propellers design is the geometrical parametrization of a single blade of the propeller.

\begin{figure}
    \centering
    \includegraphics[width=0.3\textwidth, trim={20cm 10cm 10cm 20cm}, clip]{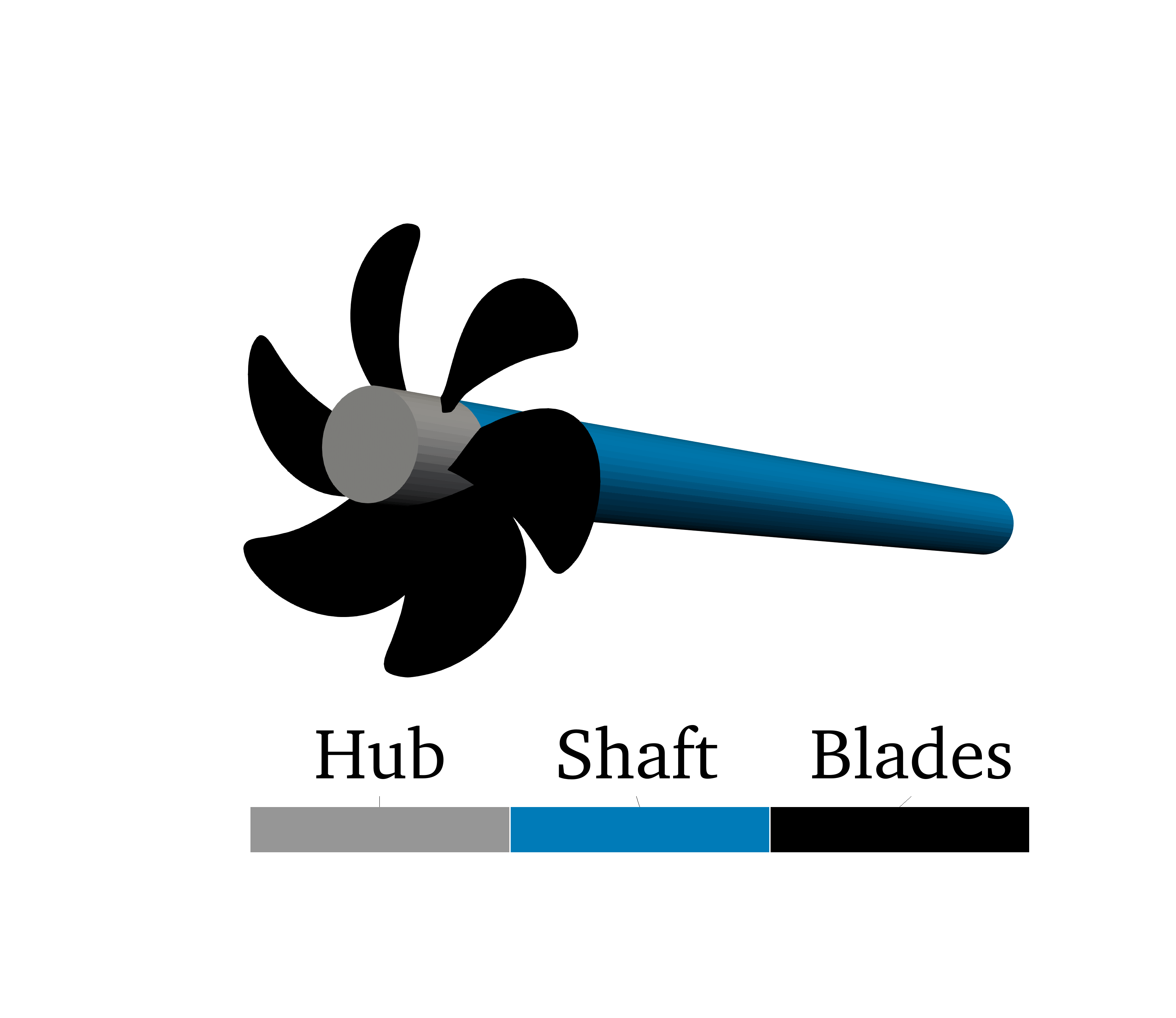}
    \caption{Generic structure of a propeller.}
    \label{subfig:ids_propeller}
\end{figure}

The generic structure of a propeller is displayed in Figure \ref{subfig:ids_propeller}, where the main components of the propeller are distinguished, the blades, the hub and the shaft.
As can be seen from Figure \ref{subfig:ids_propeller}, in this work we focus on a propeller a fixed number of blades, namely $6$ blades.

\begin{figure}[htbp]
    \centering
    \includegraphics[width=0.75\textwidth]{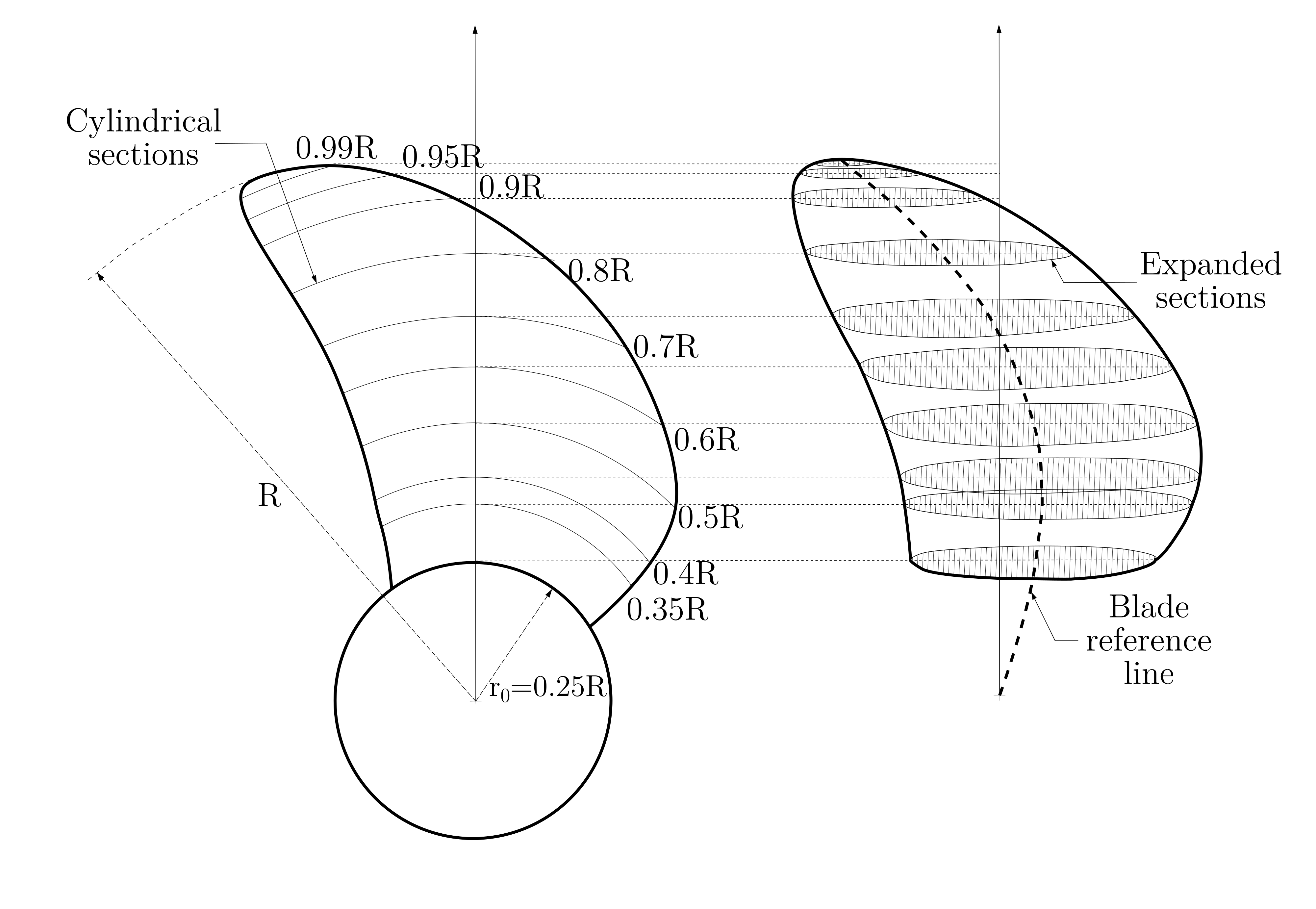}
    \caption{Radial view of a propeller's blade. The figure on the left shows the original cylinder sections of a blade, where $R$ is the propeller radius and $r_0$ is the hub radius. The figure on the right represents the projections of the corresponding sections on a flat plane.}
    \label{fig:params}
\end{figure}

\begin{figure}[htbp]
    \centering
    \subfloat[]{\includegraphics[width=0.44\textwidth]{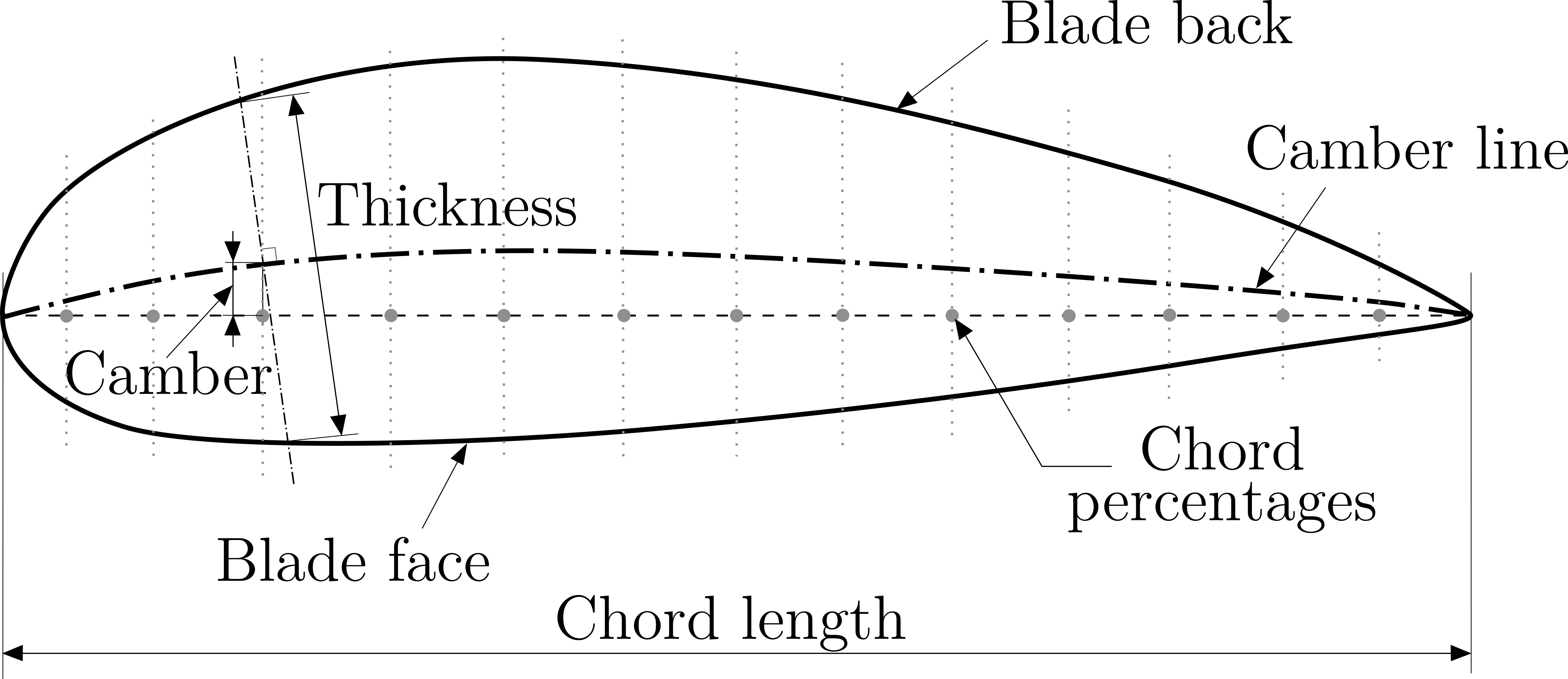} \label{subfig:section_params}}
    \subfloat[]{\includegraphics[width=0.54\textwidth]{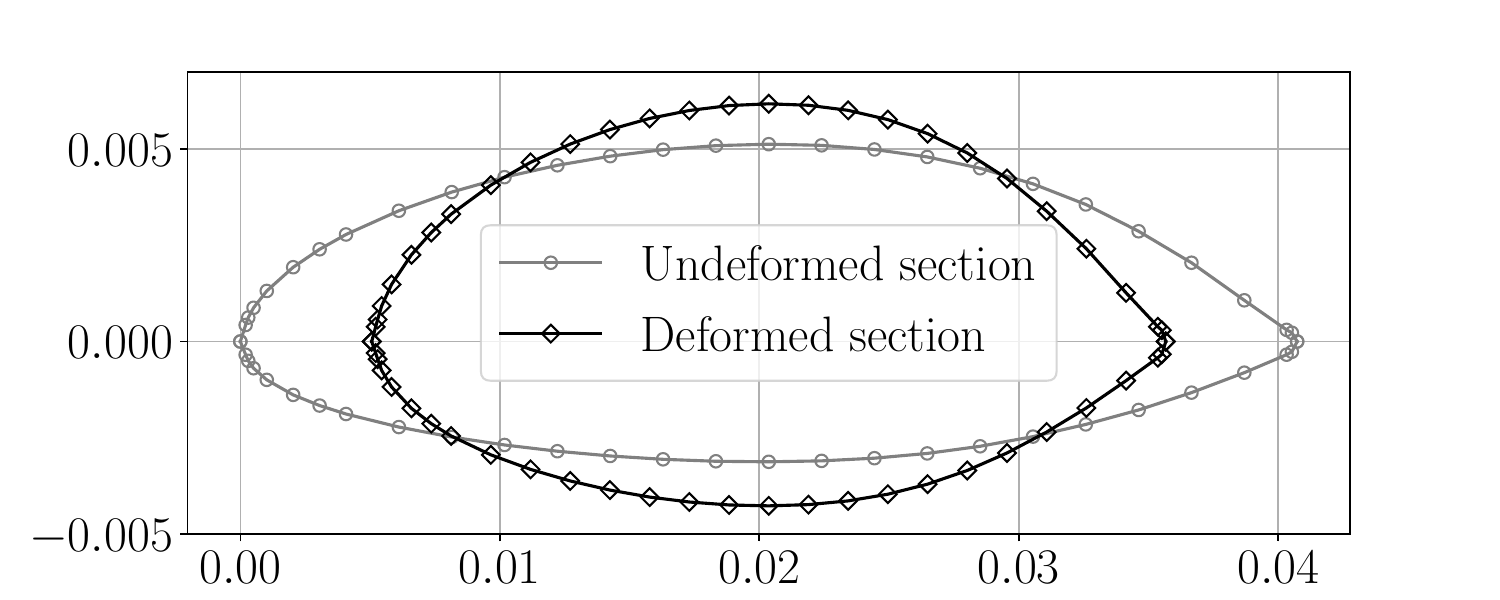} \label{subfig:section_def}}
    \caption{\protect \subref{subfig:section_params}: Graphical definition of the parameters of a blade section: chord length, thickness, and camber corresponding to the preset chord percentages.\\
    \protect \subref{subfig:section_def}: Example of deformed blade root section.}
    \label{fig:section}
\end{figure}

The blades of marine propellers are characterized by different geometrical features, which are measured on different cylindrical sections and are associated with the values of the radius at which sections are taken.
Figure \ref{fig:params} represents the values of the radii taken into account in our particular case study for blade parametrization. In particular, $r_0$ is the hub radius, whereas $R$ is the maximum radius of the propeller. As a usual approach for the study of marine propellers, the reference radii $r$ are all measured with respect to $\frac{r}{R}$, as can be seen from Figure \ref{fig:params}. In particular, the blade is characterized by four faces: the sections at radius $r_0$ and $R$, named \emph{root} and \emph{tip} respectively, the blade \emph{face} (or \emph{pressure side}) and \emph{back} (or \emph{suction side}), whose borders are indicated in Figure \ref{fig:section}. 
We can distinguish two classes of parameters: the \emph{global} parameters, associated with each radius value, i.e. pitch, rake, skew, chord length; the \emph{section} parameters, i.e. quantities defined at different chord percentages for each section, such as thickness and camber\footnote{The generality of the pipeline allows to select different and more parameters}.

The parameters which are taken into account for the blade deformation presented in this work are pitch, chord length, thickness, and camber, whereas the other parameters are fixed for all the deformed shapes.
More in detail, the pitch of a propeller is the displacement that it makes in a complete spin of $360$ degrees, we are considering a propeller with a different pitch distribution for each section considered.
The other parameters considered, i.e. chord length, thickness, and camber, are graphically defined in Figure \ref{subfig:section_params}. In particular, for each section thickness and camber are defined as lists of values associated with predefined chord percentages, as represented in Figure \ref{subfig:section_params}.
The parameters of the initial shape are reserved data\footnote{for Non-Disclosure Agreements}, but we highlight that the entire pipeline is reproducible with all the propellers defined by the above-mentioned parametrization.

\textbf{Remark} \\ It is worth remarking that we adopt the American convention in the definition of the thickness since it is measured along lines orthogonal to the camber line. In the British convention, instead, it is measured along vertical lines, orthogonal to the chord line at the specific chord percentages considered.
%The focus of the parametrization step is the extraction of the geometrical parameters starting from the cartesian coordinates of each section of the blade.

\subsection{Blade deformation}
\label{subsec:blade-def}

Once the geometrical parameters \RB{of interest} are selected, changes in the blade geometry can be easily made by modifying the parameters. In this work, 
the deformation parameters are the multiplicative factor we impose on the original geometrical parameters of the original blade. For instance, if we consider a deformation parameter equal to $\num{1.3}$ for the pitch, the resulting deformed blade will have the pitch geometrical parameter $30\%$ greater than the original blade.

The deformation rates' values are set within the intervals:
\begin{itemize}
\item $\begin{bmatrix} 0.9,&1.1 \end{bmatrix}$ (for pitch);
\item $\begin{bmatrix} 0.8,&1.2 \end{bmatrix}$ (for camber);
\item $\begin{bmatrix} 0.7,&1.3 \end{bmatrix}$ (for chord length);
\item $\begin{bmatrix} 0.7,&1.3 \end{bmatrix}$ (for thickness).
\end{itemize}
\RA{The ranges have been chosen such that the parametric blades satisfy the structural feasibility constraint.}
Figure \ref{subfig:section_def} displays an example of deformation of the root section of the blade when the deformation parameter is $\mu = \begin{bmatrix}
1&0.95&0.75&1.27
\end{bmatrix}$, where the components refer to pitch, camber, chord length, and thickness, respectively.
Such deformation is obtained by employing the Python package \verb+BladeX+ \cite{gadalla19bladex, bladex}. The software builds the blades by creating a Non-Uniform Rational Basis Spline (NURBS) surface~\cite{piegl1995tessellating, piegl1996nurbs} passing through all the (deformed) blade sections.
In this representation, the shape of the surfaces is defined by the control points belonging to the UV plane, which enables the following steps of the pipeline thanks to the mapping between such a (2D) plane and the 3D final surface.

\textbf{Remark} \\
\RB{It is worth highlighting that we considered the same deformation rates for each section and each chord percentage, i..e,} we are \RB{not} taking into \RB{account} different deformations for different sections and/or chord coordinates. Therefore, we are considering $\nparams{}=4$ parameters for the whole blade.\\

\begin{figure}[htbp]
    \centering
    \subfloat[]{\includegraphics[width=0.3\textwidth, trim={5cm 0 5cm 0}, clip]{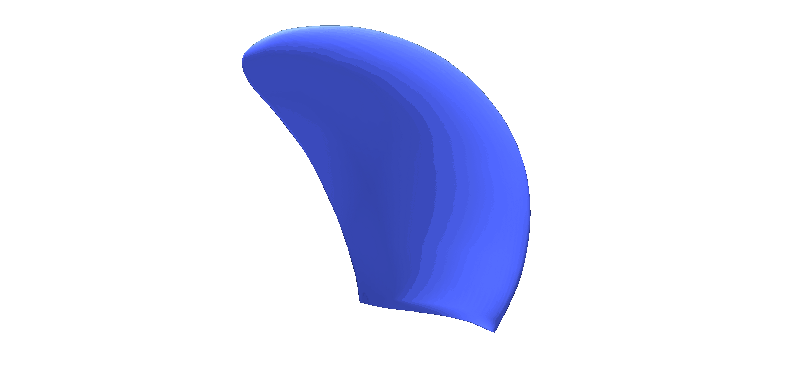}}\subfloat[]{\includegraphics[width=0.3\textwidth, trim={5cm 0 5cm 0}, clip]{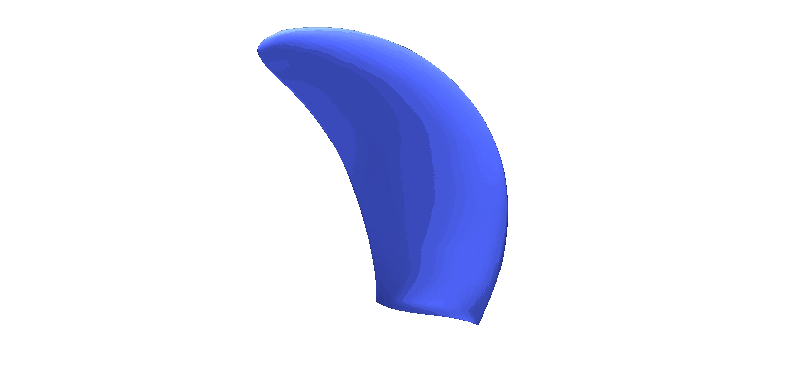}}\subfloat[]{\includegraphics[width=0.3\textwidth, trim={5cm 0 5cm 0}, clip]{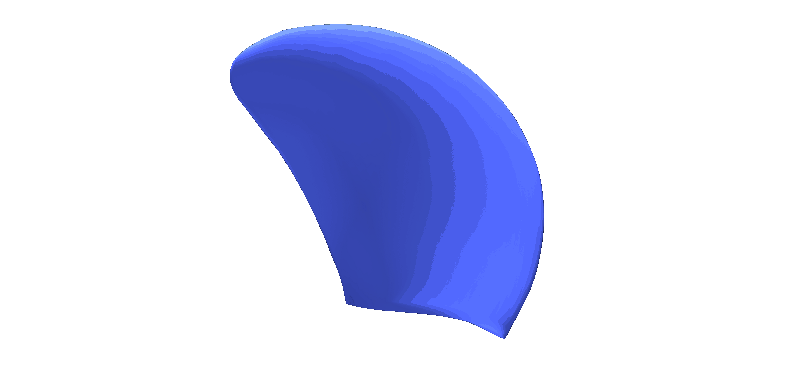}}
    \caption{Graphical view of three deformed blades.}
    \label{fig:blades-def}
\end{figure}

In the project here presented, we considered $\nsnaps{}=216$ deformed blades, where $200$ blades are obtained imposing deformations with parameters selected with uniform random distribution in the intervals defined above, and $16$ blades are deformed considering all the combinations at the extremes of the intervals.
Three examples of deformed shapes are displayed in Figure \ref{fig:blades-def}.

\subsection{Mesh deformation}
\label{subsec:mesh-def}

A challenging step in the pipeline of this project is to perform the mesh deformation.
Indeed, after a large number of deformed blades is obtained, as described in Section \ref{subsec:param}, the computational mesh (schematized in Figure \ref{fig:domain}) has to be deformed as well, according to each blade deformation.
However, it is important to highlight one important requirement that all the deformed meshes have to satisfy: the mesh topology, as well as the number of degrees of freedom, should be preserved in all the \fom{} simulations.
The statement above is a pre-requisite for the construction of the reduced order model since all the pressure and velocity high-fidelity solutions should have the same degrees of freedom.

A mesh that has the same number of points and the same topology as the \fom{} simulation with the undeformed blades cannot be obtained from the \emph{OpenFOAM} tool of mesh generation. Therefore, the mesh deformation is performed in this work by exploiting an interpolation technique. In particular, we followed two fundamental steps to obtain the points of each deformed mesh: the deformation of the propeller's shaft and the deformation of the entire mesh, described in detail in the following parts. The images reported in Subsections \ref{subsec:shaft-def}, \ref{subsec:def-all-mesh} and \ref{subsec:results-def} display results for the deformed blade with parameter $\mu = \begin{bmatrix}
    1&0.72&0.75&0.95
\end{bmatrix}$.

Before specifying the operative details of the deformation, we recall in the following paragraph the basic concepts of the Radial Basis Function (\rbf{}) interpolation technique when used for mesh deformation.

\paragraph{\rbf{} for mesh deformation}
The general problem is to find the deformed counterpart $\mathbf{x}^{\text{def}}$ of the original undeformed mesh $\mathbf{x}^{\text{undef}}$. As starting point we should consider a certain number $N_{\text{control}}$ of undeformed and deformed control points $\mathbf{c}^{\text{undef}}$ and $\mathbf{c}^{\text{def}}$. The deformed mesh is found with the following expression:
\begin{equation}
x_j^{\text{def}} = f(x_j^{\text{undef}}) = \sum_{i=1}^{N_{\text{control}}} \alpha_i \phi(\| x_j^{\text{undef}} - c_i^{\text{undef}}\|), \, j=1,\dots, N_{\text{points}}\,,
    \label{eq:apply-rbf}
\end{equation}
where $\phi_i = \phi(\| x_j^{\text{undef}} - c_i^{\text{undef}}\|)$ are the radial basis functions, which can have different shapes, such as thin plate splines, multiquadric or inverse multiquadric. In our case, we consider the \emph{thin plate splines}, thus: $\phi(r)=r^2\text{log}(r)$, where $r$ is the radius of the basis function, i.e. $r=\| x_j^{\text{undef}} - c_i^{\text{undef}}\|$.
\textbf{Remark} \\ We highlight here that choosing a proper \rbf{} kernel would ensure the preservation of the mesh quality and topology throughout the discretized domain.
The weights $\alpha_i$ are found by \emph{training} the interpolation with some known control points, that have to satisfy the following conditions:
\begin{equation}
c_j^{\text{def}} = f(c_j^{\text{undef}}) = \sum_{i=1}^{N_{\text{control}}} \alpha_i \phi(\| c_j^{\text{undef}} - c_i^{\text{undef}}\|), \, j=1,\dots, N_{\text{control}}\,.
    \label{eq:train-rbf}
\end{equation}

\subsubsection{Shaft deformation}
\label{subsec:shaft-def}
\begin{figure}[htbp]
    \centering
    \subfloat[]{\includegraphics[width=0.3\textwidth]{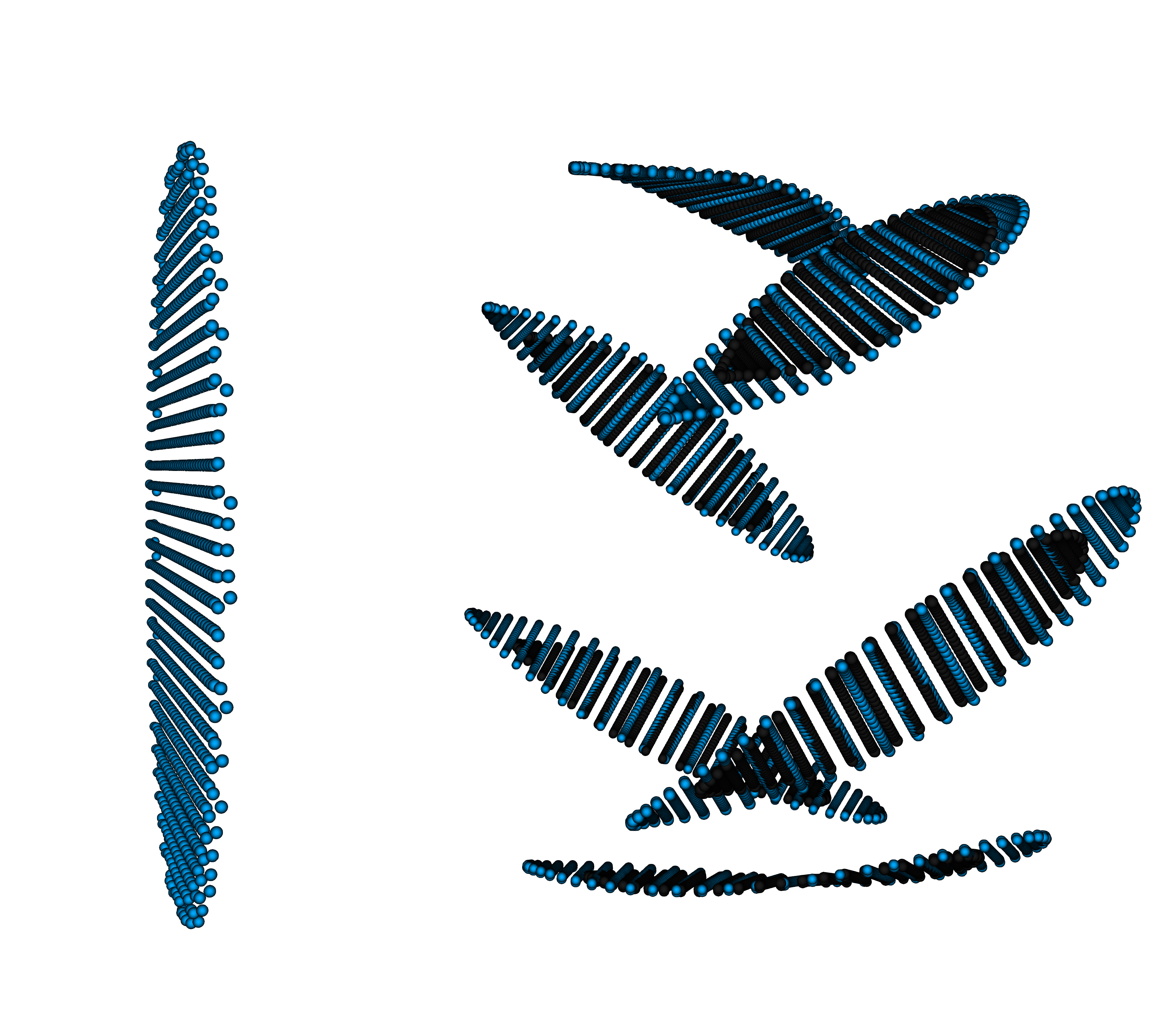}\label{subfig:control-points-shaft}}
    \subfloat[]{\includegraphics[width=0.65\textwidth]{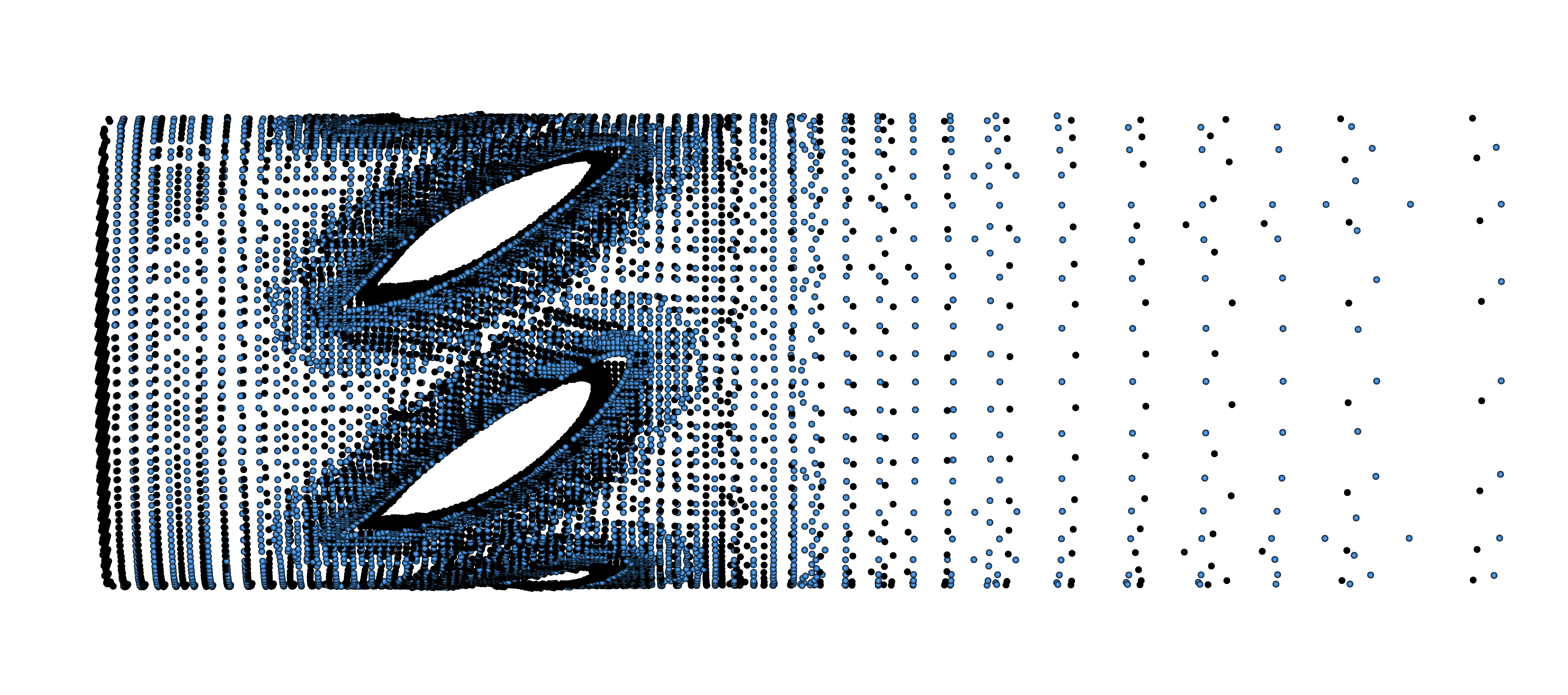}\label{subfig:all-points-shaft}}
    \caption{\protect \subref{subfig:control-points-shaft}: Zoom on the control points on the propeller's shaft. \protect \subref{subfig:all-points-shaft} Shaft points. Undeformed and deformed shaft points are represented in blue and black, respectively.} 
    \label{fig:shaft}
\end{figure}

For each deformed blade, the following list of passages is followed to obtain the deformed shaft:
\begin{enumerate}
    \item identification of $N_{\text{shaft}}$ points on the root face of the undeformed and deformed blades. Being the surface defined by NURBS, the points are generated on the reference UV plane and mapped to the blade surface;
    \item collection of the shaft points which remain unchanged in deformation, i.e. the two bases of the shaft. The union of the undeformed shaft bases and of the undeformed blades' root provides the \emph{undeformed shaft control points}, whereas the union of the undeformed shaft bases and of the deformed blades' root provides the \emph{deformed shaft control points};
    \item training of a \emph{Radial Basis Function} (\rbf{}) interpolating technique with the \emph{undeformed} and \emph{deformed shaft control points} (as in \eqref{eq:train-rbf}), and then extraction of the deformed shaft by applying the \rbf{} to the undeformed points of the shaft lateral surface (as in \eqref{eq:apply-rbf}).
    The control points and the final mesh points on the shaft are represented in Figure~\ref{fig:shaft}.
\end{enumerate}

To conclude, we highlight that such an operation is needed in order to allow the correct deformation of the whole mesh. \RB{If the points on the shaft surface are fixed, this would lead} indeed to unexpected behavior when the RBF system is assembled over all the surface points. \RB{In fact, in this case} the blade points are moved, but not the adjacent ones belonging to the shaft. The propedeutic procedure here described induces indeed a deformation on the shaft points, translating them on the cylindrical surface accordingly to the blade root deformation.

\subsubsection{Global mesh deformation}
\label{subsec:def-all-mesh}
\begin{figure}[htbp]
    \centering
    \subfloat[]{\includegraphics[width=0.5\textwidth]{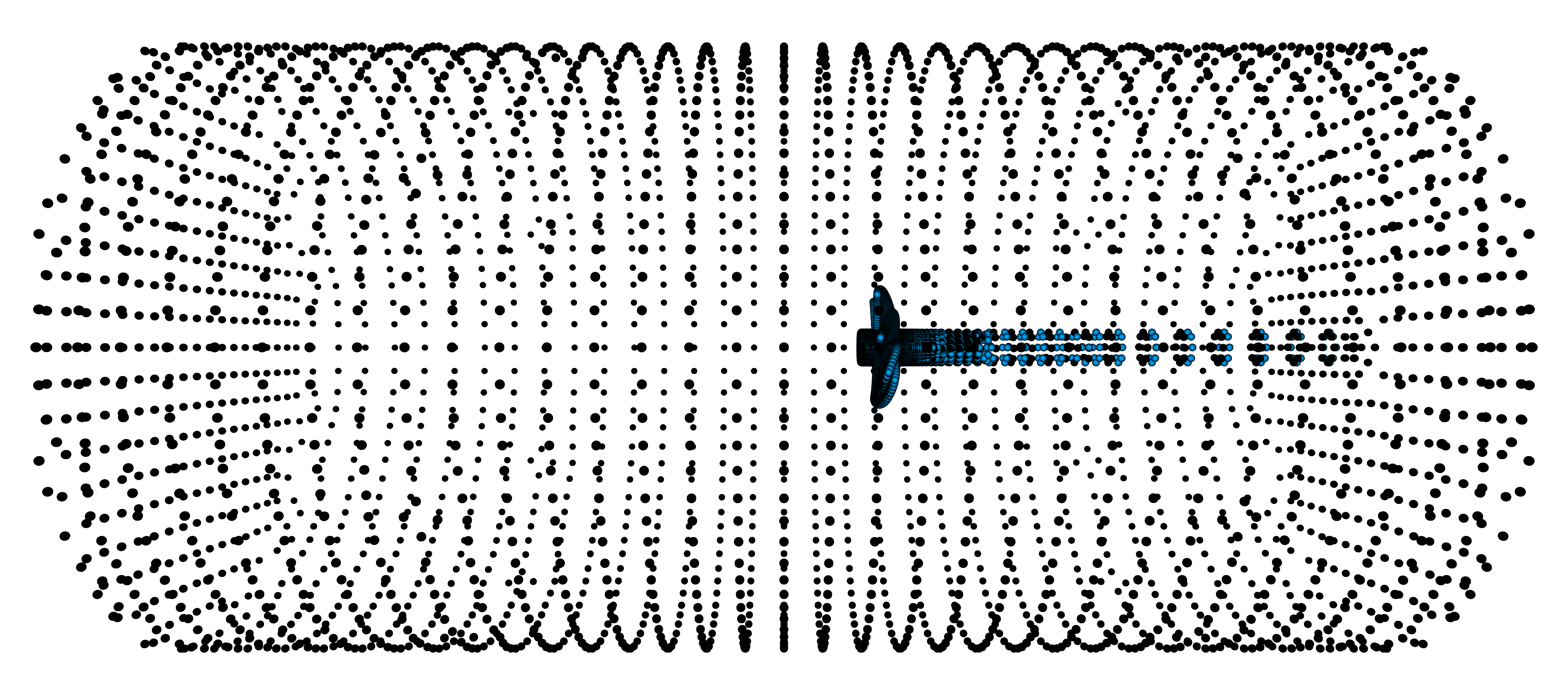}\label{subfig:all_control_points}}
    \subfloat[]{\includegraphics[width=0.45\textwidth]{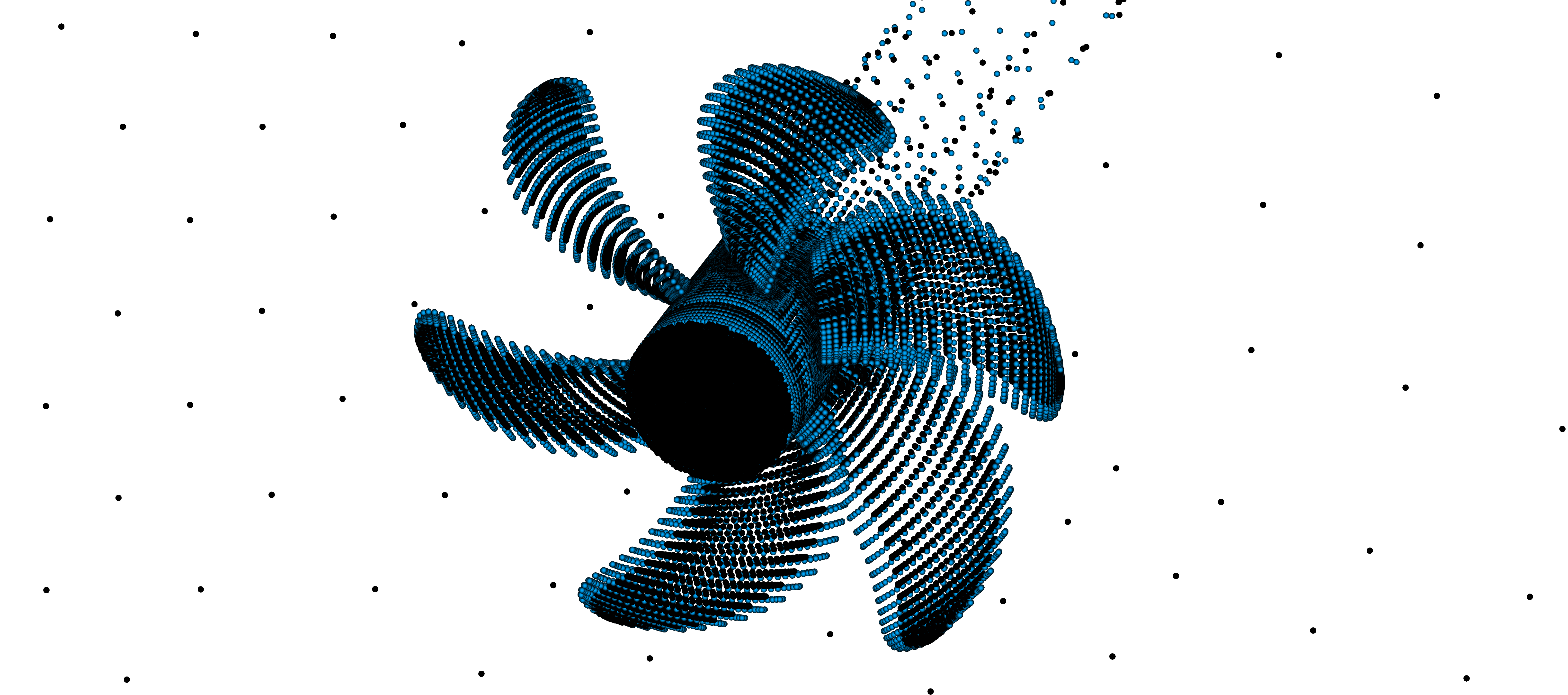}\label{subfig:zoom_blades}}
    \caption{\protect\subref{subfig:all_control_points}: Representation of all the control points in the mesh, including points at boundaries and on the propeller's surface; \protect\subref{subfig:zoom_blades}: Zoom on the undeformed (blue) and deformed (black) control points on blades.} 
    \label{fig:deformation_mesh}
\end{figure}

The deformation of the whole mesh is obtained with the following steps:
\begin{enumerate}
    \item identification of $N_{\text{quadrature}}$ Gauss quadrature nodes on the faces of both the undeformed and deformed blades; we name these nodes as \emph{undeformed} and \emph{deformed blades' control points}. Here we chose as control points the Gauss quadrature nodes on the blades, but other selections that ensure a proper geometrical characterization of the blades are possible.
    \item collection of all the points of the undeformed mesh which belong to the boundaries, i.e. the lateral surface of the outer cylinder $\Gamma_{\text{outer}}$, the inlet $\Gamma_{\text{inlet}}$, the outlet $\Gamma_{\text{outlet}}$, and the surfaces of the propeller shaft $\Gamma_{\text{shaft}}$. The union of the undeformed blades' control points and the undeformed boundaries will provide the \emph{undeformed mesh control points}, while the union of the deformed blades' control points and the undeformed boundaries will provide the \emph{deformed mesh control points};
    \item training of a second \rbf{} interpolation with the undeformed and deformed mesh control points (as reported in \eqref{eq:train-rbf}); then, the \rbf{} is applied on the entire internal mesh (again, as in \eqref{eq:apply-rbf}) to obtain the deformed internal mesh, which is represented in Figure \ref{fig:deformation_mesh}.
\end{enumerate}

\textbf{Remark} \\ It is important to remark that, if we want to perform a deformation of a certain domain, there is the need to include in the control points the boundaries of the domain, i.e. the bases of the shaft in point (1) and the boundaries of the whole domain in point (2).

\subsubsection{Results for mesh deformation}
\label{subsec:results-def}
This part is dedicated to the graphic results of the mesh deformation step. In particular, here we show different mesh views for the original and the deformed blades and for the meshes of the corresponding \fom{} simulations.
For a detailed description of the \fom{} setting we refer the reader to Section \ref{sec:fom}.

\begin{figure}[htbp]
    \centering
    \subfloat[Mesh on undeformed blade surface]{\includegraphics[width=0.35\textwidth]{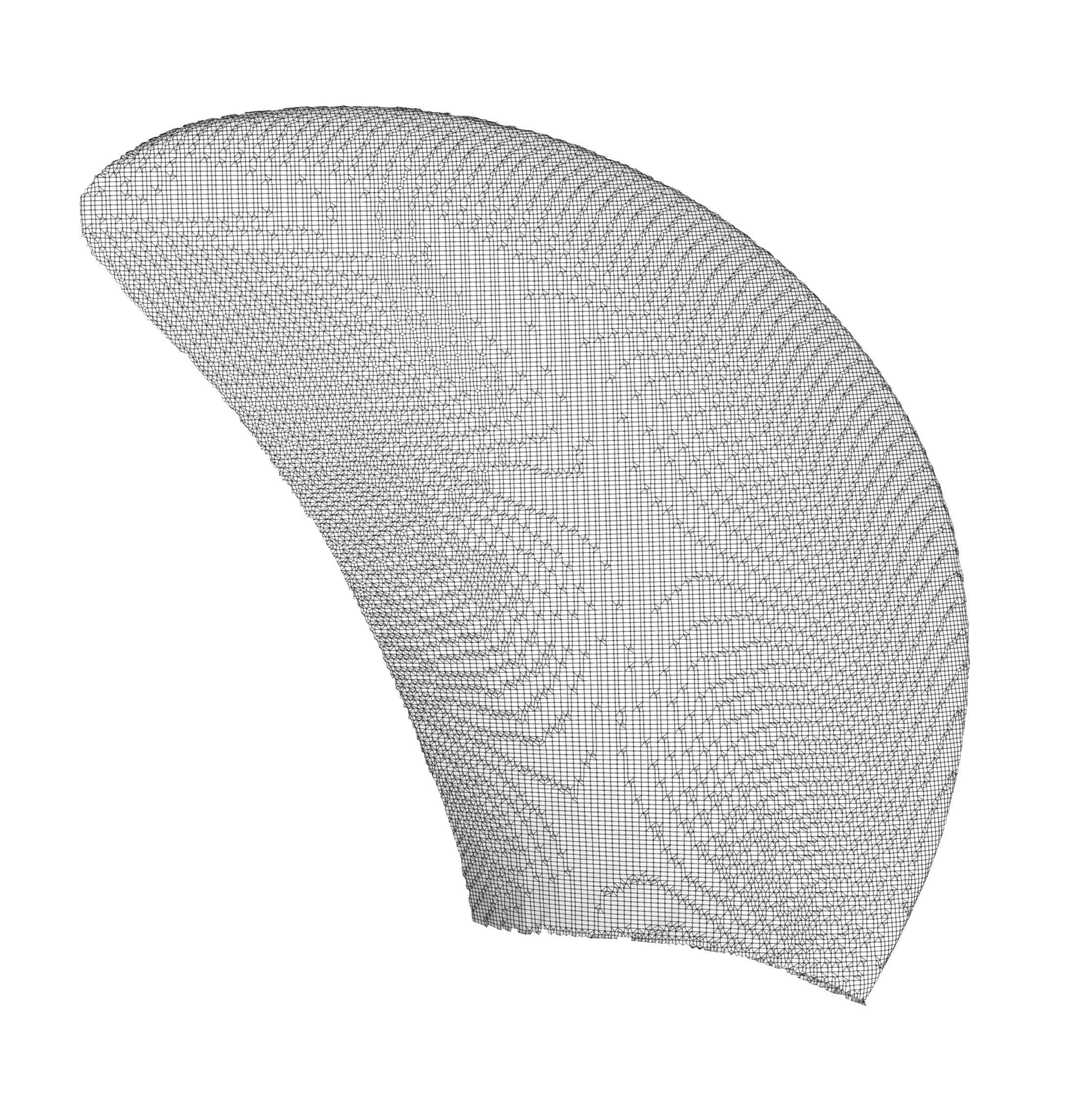}}
    \subfloat[Mesh on deformed blade surface]{\includegraphics[width=0.35\textwidth]{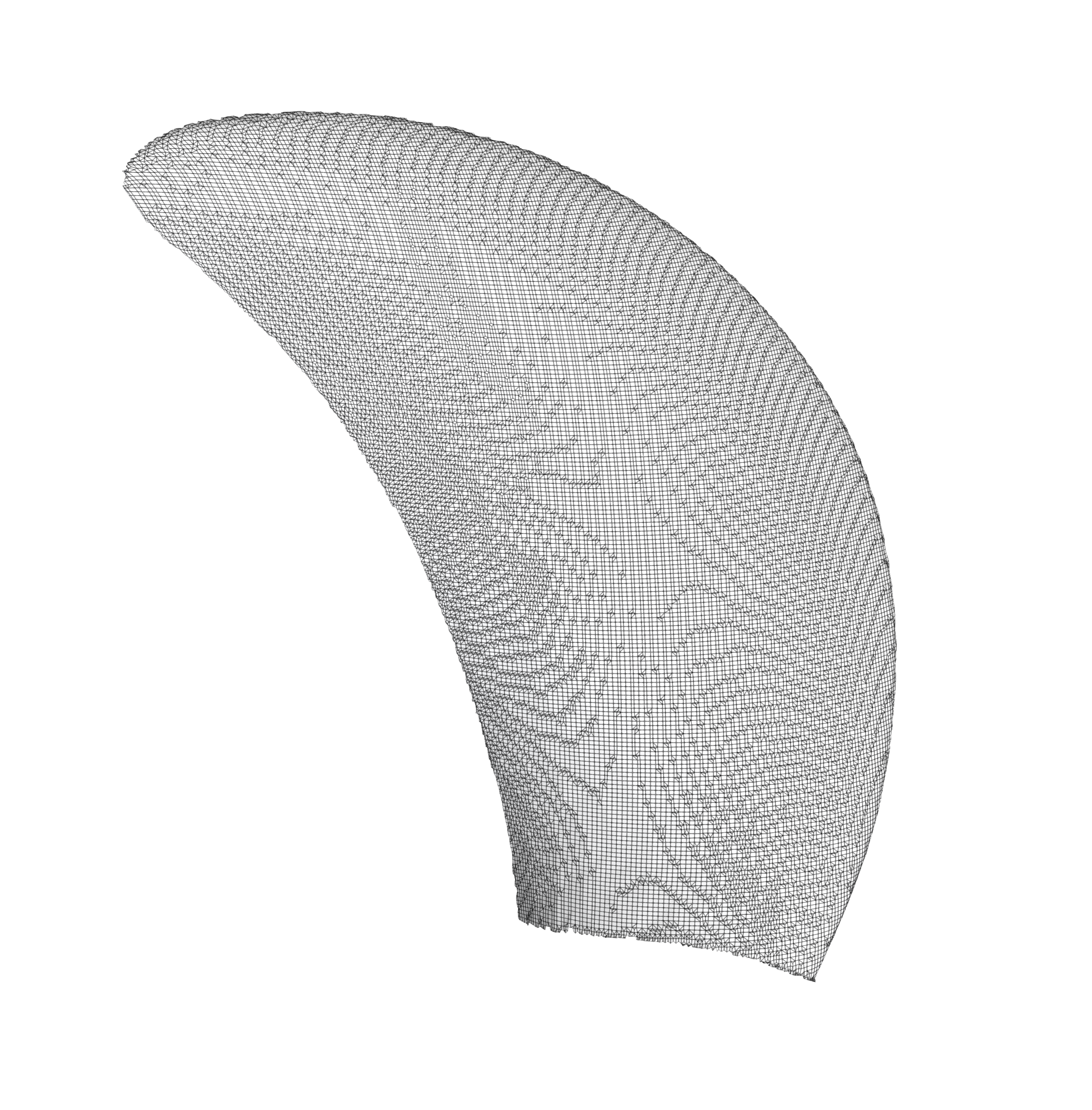}}
    \caption{Mesh on the surface of the original undeformed blade \RB{and of a selected deformed blade}.}
    \label{fig:mesh-blade}
\end{figure}

\begin{figure}[htbp]
    \centering
    \subfloat[Undeformed mesh slice]{\includegraphics[width=0.35\textwidth]{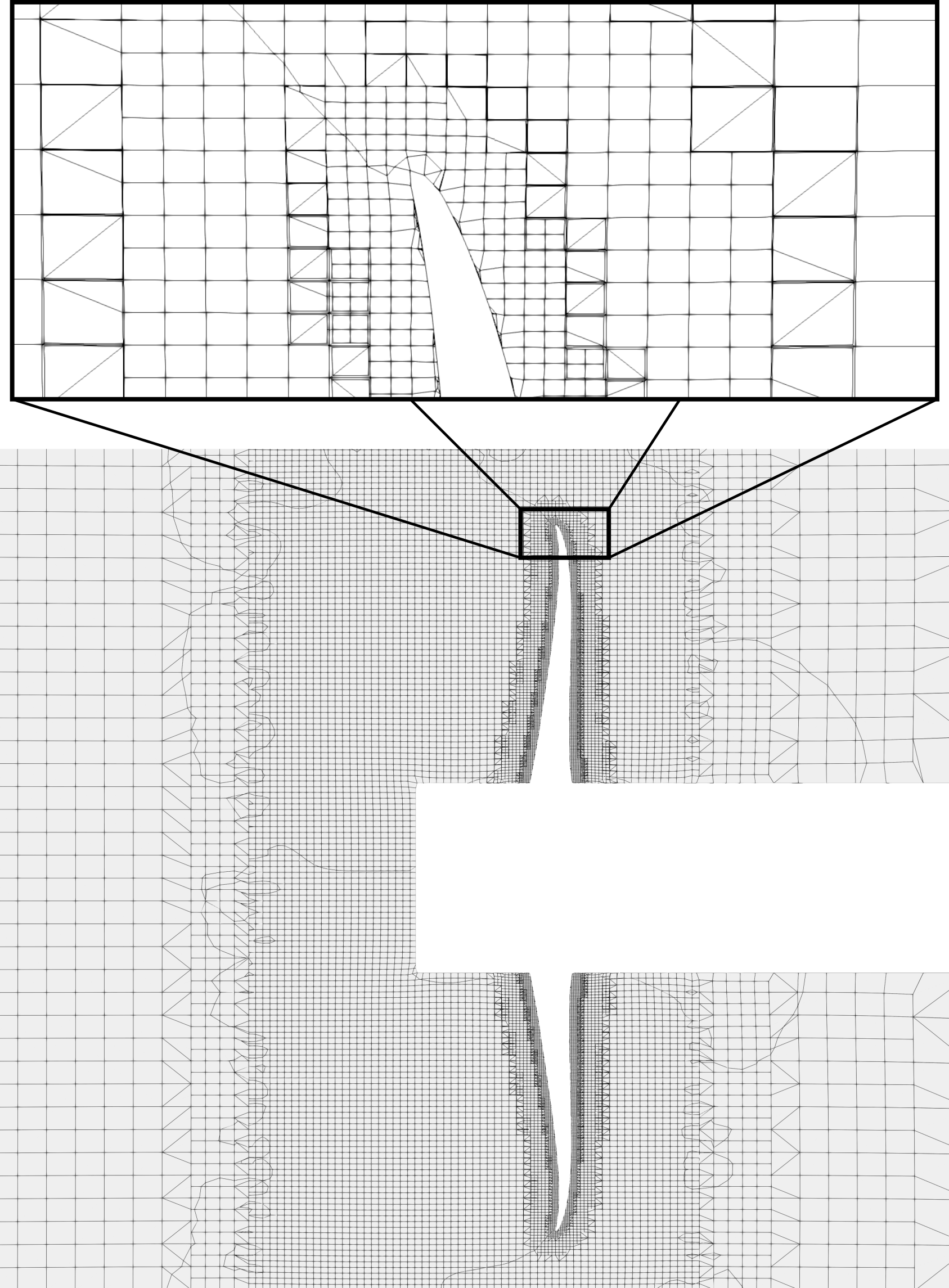}}
    \quad
    \subfloat[Deformed mesh slice]{\includegraphics[width=0.35\textwidth]{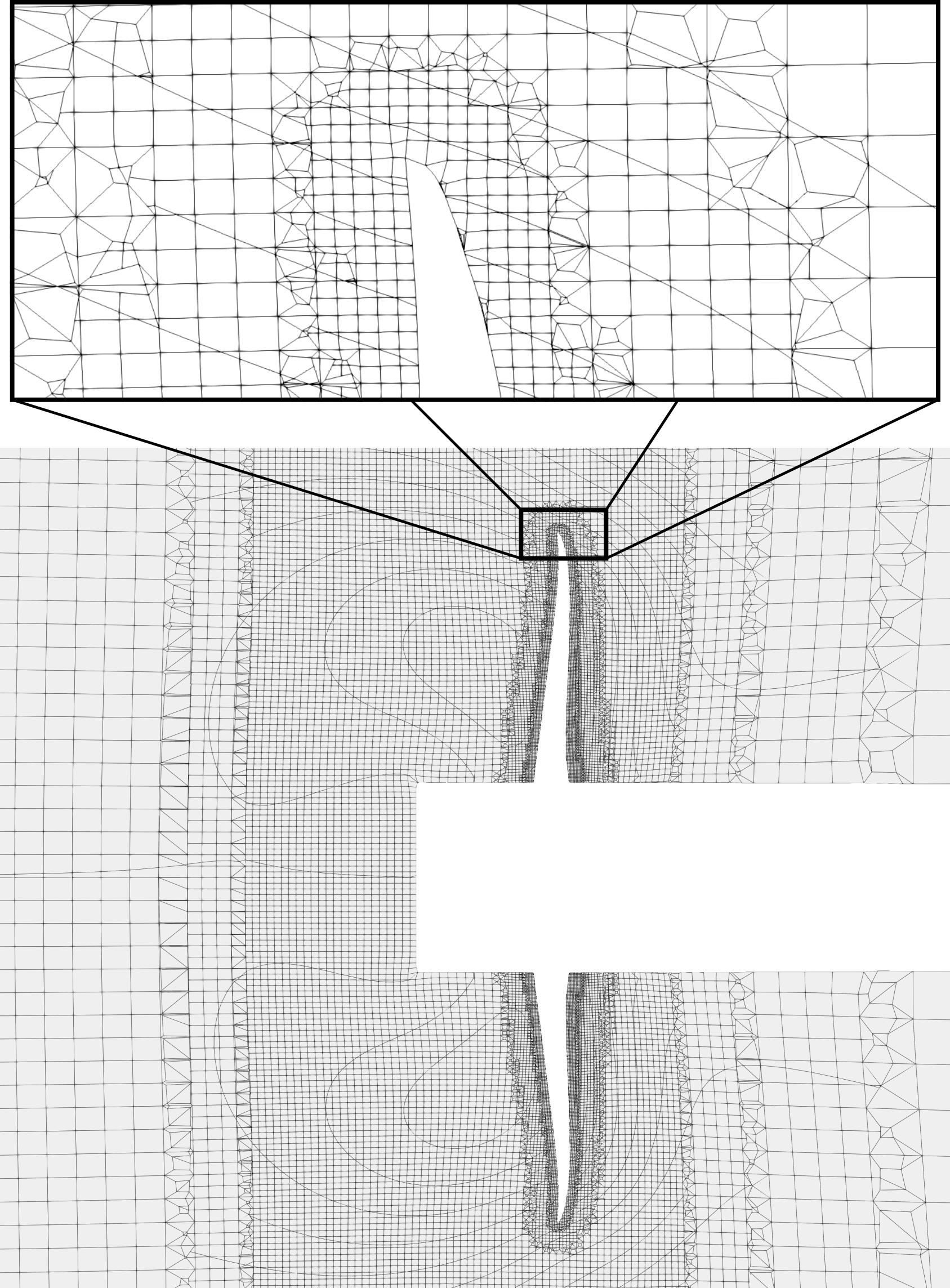}}
    \caption{Zoom of a slice of the mesh with corresponding normal parallel to the propeller axis.}
    \label{fig:mesh-plane1}
\end{figure}

\begin{figure}[htbp]
    \centering
    \subfloat[Undeformed mesh slice]{\includegraphics[width=0.35 \textwidth]{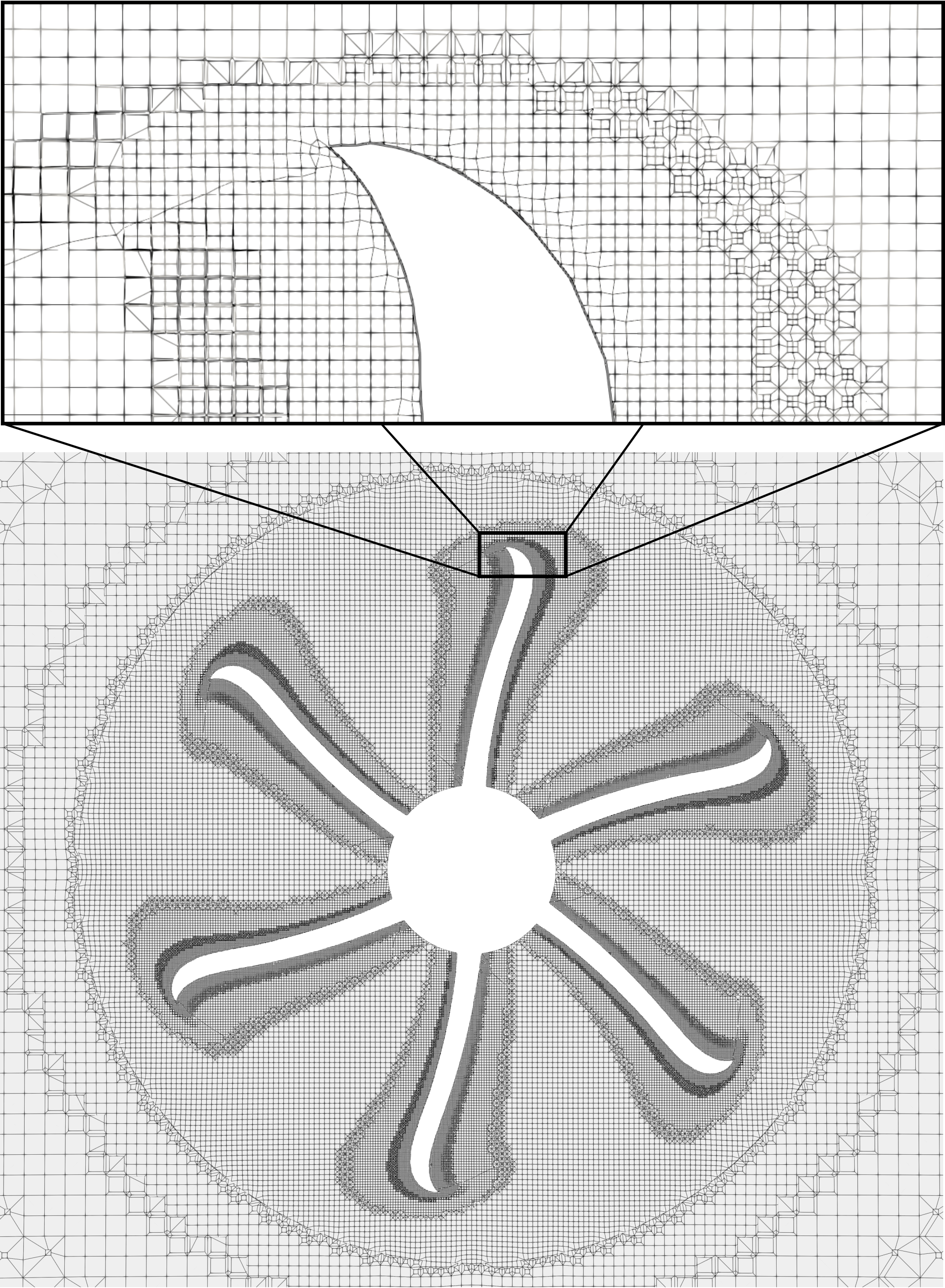}}
    \quad
    \subfloat[Deformed mesh slice]{\includegraphics[width=0.35 \textwidth]{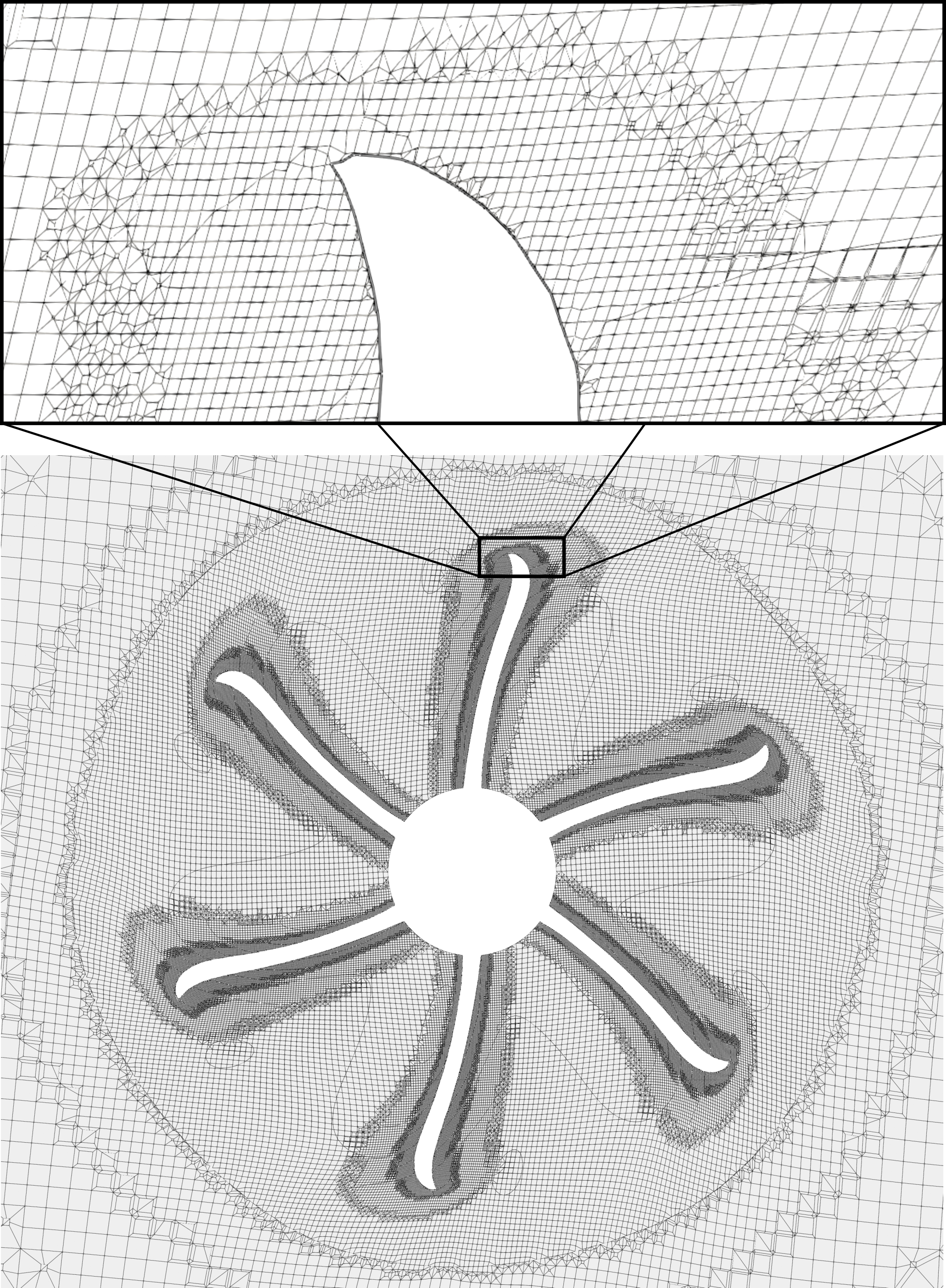}}
    \caption{Zoom of a slice of the mesh with corresponding normal orthogonal to the propeller axis.}
    \label{fig:mesh-plane2}
\end{figure}

Figure~\ref{fig:mesh-blade} displays the mesh on the surfaces of the undeformed and deformed blades, Figures \ref{fig:mesh-plane1} and~\ref{fig:mesh-plane2} represent the slices of the mesh on two different planes, parallel and orthogonal to the propeller axis, respectively.
All the deformed meshes obtained through the above-mentioned interpolation technique provide good results in terms of mesh quality.

The time needed to perform the mesh deformation using a \rbf{} interpolation method is $\sim 2-2.5$ hours in a serial computation~\footnote{The mesh deformation is performed using one processor core only on SISSA HPC cluster Ulysses (200 TFLOPS, 2TB RAM, 7000 cores).}.

\newpage

\section{Full Order Model}
\label{sec:fom}
A necessary step in our pipeline is the setting of the Full Order Model (\fom{}) for high-fidelity simulations. \RA{The mathematical model here considered is the Unsteady Reynolds Averaged Navier--Stokes (U-RANS) Equations, that are numerically discretized and solved making use of} the open-source software \emph{OpenFOAM}~\cite{ofsite}, which exploits the finite-volume method~\cite{moukalled2016finite}. \RA{Subsection \ref{subsec:rans} describes the U-RANS approach coupled with the turbulence modeling, whereas Subsection \ref{subsec:mrf} explains the technique used to generate the rotation of the propeller in the FOM simulations.} 

\subsection{The U-RANS approach and the turbulence modeling}
\label{subsec:rans}

In this paper, we adopt the following notation: $\Omega \in \mathbb{R}^d$ with
$d=2$ or $3$ is the fluid domain, $\Gamma$ its boundary, $t \in [0,T]$ the time, $\mathbf{u}=\mathbf{u}(\mathbf{x},t)$ the velocity field, $p=p(\mathbf{x},t)$ the normalized pressure scalar field divided by the fluid density, and $\nu$ the fluid kinematic viscosity. 

\begin{figure}[h!]
    \centering
    \includegraphics[width=0.65\textwidth]{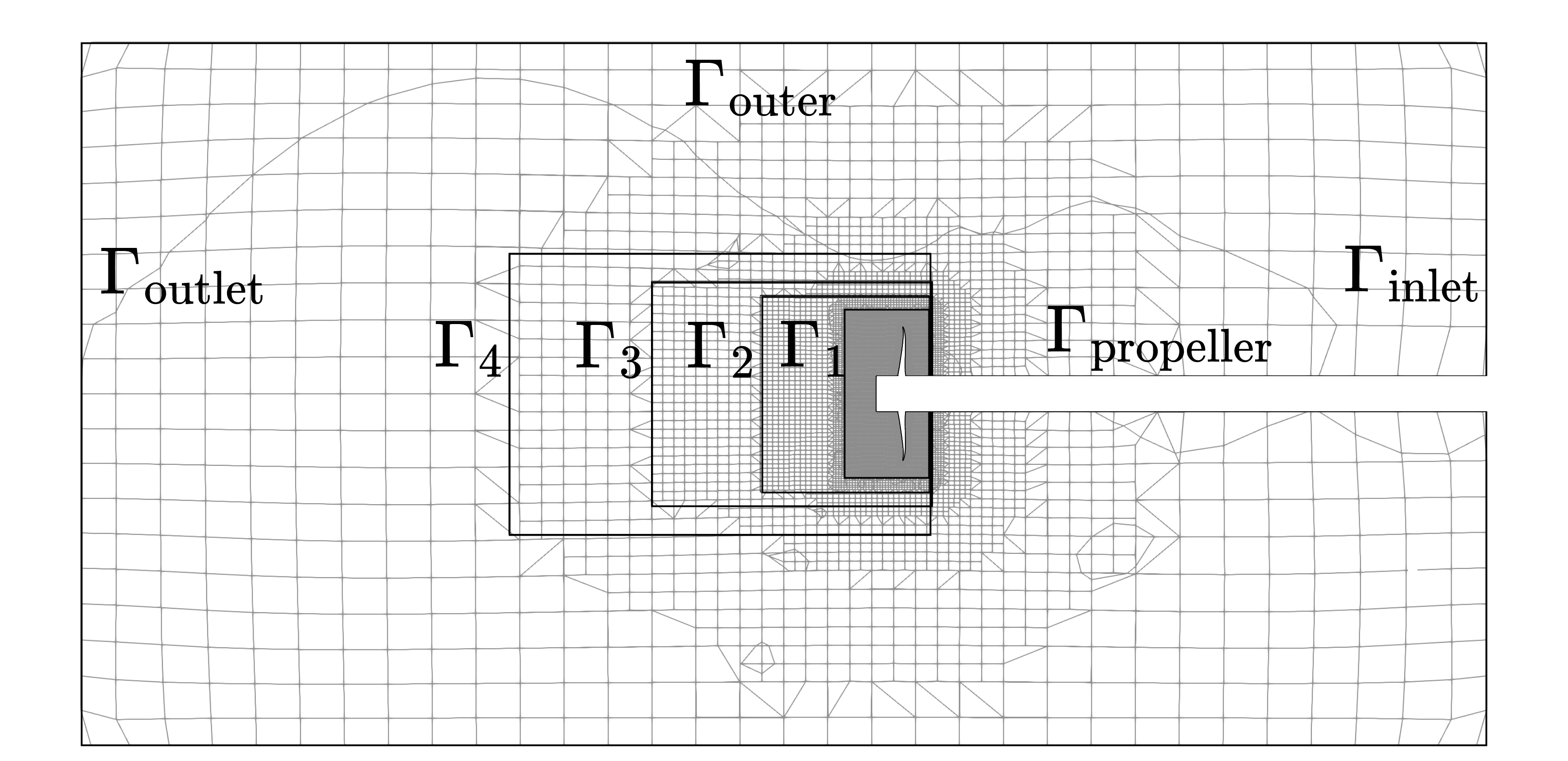}
    \caption{A slice of the computational domain $\Omega$ and mesh of the \fom{} in \emph{OpenFOAM}}
    \label{fig:domain}
\end{figure}

Figure \ref{fig:domain} represents a slice of our computational domain and mesh. In particular, the mesh is built inside the outer cylinder, whose bases are the inlet $\Gamma_{\text{inlet}}$, the outlet $\Gamma_{\text{outlet}}$ and the lateral surface $\Gamma_{\text{outer}}$.
The surface of the propeller is named here $\Gamma_{\text{propeller}}$. In the internal domain, four different cylinders are built, whose surfaces are indicated as $(\Gamma_i)_{i=1}^4$.
We introduce here the internal cylinders in order to generate different mesh refinements; in particular, the mesh gets gradually more refined from the outer surface to the blades' surface in order to provide an accurate reconstruction of the flow fields acting on the blades, as can be seen from Figure \ref{fig:domain}.

\begin{comment}
The incompressible \nse{} in the described domain is written as follows:

\begin{subnumcases}{\label{NSE}}
\frac{\partial \mathbf{u}}{\partial t}=-\nabla \cdot (\mathbf{u} \otimes \mathbf{u})+\nabla \cdot \nu \left(\nabla \mathbf{u} + (\nabla \mathbf{u})^T \right)-\nabla p & in $\Omega \times [0,T]\, , $ \label{mom_NSE}\\
\nabla \cdot \mathbf{u} = \mathbf{0} & in $\Omega \times [0,T]\, ,$ \label{cont_NSE}\\
+ \text{ boundary conditions }& on  $\Gamma \times [0,T]\, ,$ \label{bound_NSE}\\
+ \text{ initial conditions }& in  $(\Omega,0)\, .$ \label{init_NSE}
\end{subnumcases}
\end{comment}

\RA{We briefly recall the basic concepts of the U-RANS approach, that is used in this work for the FOM simulations.
The main hypothesis that characterizes the RANS approach is the \emph{Reynolds decomposition}~\cite{reynolds1895iv}.
This theory is based on the assumption that each flow field can be expressed as the sum of its mean and fluctuating parts. Considering a generic space and time-dependent field $\sigma(\mathbf{x}, t)$:
\[
\sigma(\mathbf{x}, t) = \overline{\sigma}(\mathbf{x}, t) + \sigma^{\prime}(\mathbf{x}, t).
\]

The U-RANS formulation consists in a time-averaged version of the NSE, that can be written as follows:
\begin{equation}
\begin{cases}
\dfrac{\partial \overline{u_i}}{\partial x_i} = 0 ,\\
    \dfrac{\partial \overline{u_i}}{\partial t}+ \overline{u_j} \dfrac{\partial \overline{u_i} }{\partial x_j}=-\dfrac{\partial \overline{p}}{\partial x_i}+ \dfrac{\partial( 2 \nu\overline{\mathbf{E}}_{ij}- \mathcal{R}_{ij})}{\partial x_j},
    \label{URANS}
\end{cases}
\end{equation}
where the Einstein notation has been adopted, $\mathcal{R}_{ij}=\overline{u'_i u'_j}$ is the Reynolds stress tensor, and $\overline{\mathbf{E}}_{ij}=\dfrac{1}{2}\left(\frac{\partial \overline{u_i}}{\partial x_j} + \frac{\partial \overline{u_j}}{\partial x_i}\right)$ is the averaged strain rate tensor.

The U-RANS formulation in \eqref{URANS} needs to be coupled with a turbulence model in order to close system \eqref{URANS}. In particular, here we adopt the $\kappa-\omega$ \emph{Shear Stress Transport} (SST) model~\cite{menter1994two}.
}
It belongs to the class of \emph{eddy viscosity models}, whose main assumption is the Boussinesq hypothesis, i.e. that the turbulent stresses are related to the mean velocity gradients \RA{as follows:
\[
 \mathcal{R}_{ij}=2 \nu_t \mathbf{E}_{ij} - \dfrac{2}{3} \kappa \delta_{ij},
    \label{bouss}
\]
where $\kappa=\frac{1}{2}\overline{u'_i u'_i}$ is the turbulent kinetic energy, $\nu_t$ is the eddy viscosity.

The final formulation is expressed as follows:

\begin{subnumcases}{\label{U-RANS-turb}}
    \dfrac{\partial \overline{\mathbf{u}}}{\partial t}+ \nabla \cdot (\overline{\mathbf{u}} \otimes \overline{\mathbf{u}}) = \nabla\cdot\left[ -\overline{p} \mathbf{I}+(\nu + \nu_t)\left( \nabla\overline{\mathbf{u}} +(\nabla \overline{\mathbf{u}})^T\right)\right] &in $\Omega \times [0,T]$, \label{rans1}\\
    \nabla \cdot\overline{\mathbf{u}}=0 & in $\Omega \times [0,T]$, \label{rans2}\\
    + \text{ Boundary conditions} & on $\partial \Omega \times [0,T]$, \label{bound_NSE}\\
    + \text{ Initial conditions} & in $(\Omega,0)$ \label{init_NSE}.
\end{subnumcases}
}
The boundary and initial conditions \eqref{bound_NSE} and \eqref{init_NSE} are here reported.
For the pressure field, we consider:
\begin{equation*}
    \begin{cases}
        p(\mathbf{x}, 0)=0, \quad \mathbf{x} \in \Omega, \\
        \dfrac{\partial p}{\partial \mathbf{n}} (\mathbf{x}, t) = 0, \quad \mathbf{x} \in \Gamma_{\text{inlet}}, \Gamma_{\text{outer}},\\
        p(\mathbf{x}, t)=0, \quad \mathbf{x} \in \Gamma_{\text{outlet}}.
    \end{cases}
\end{equation*}
For the velocity field:
\begin{equation*}
    \begin{cases}
        \mathbf{u}(\mathbf{x}, 0) = \mathbf{0}, \mathbf{x} \in \Omega,\\
        \mathbf{u}(\mathbf{x}, t) = \mathbf{u}_0(\mathbf{x})=(0, u_0, 0), \quad \mathbf{x} \in \Gamma_{\text{inlet}}\\
        \dfrac{\partial\mathbf{u}}{\partial \mathbf{n}}(\mathbf{x}, t) = 0, \quad \mathbf{x} \in \Gamma_{\text{outlet}}\\
        \mathbf{u}(\mathbf{x}, t) = \mathbf{0}, \mathbf{x} \in \Gamma_{\text{outer}}, \Gamma_{\text{propeller}}.\\
    \end{cases}
\end{equation*}
In our case, the value of the inlet velocity $u_0$ is set such that the advance ratio $J=\frac{u_0}{nD}=0.85$, where $D=2R$ (if we refer to $R$ in Figure \ref{fig:params}) is the diameter of the propeller, and $n=15$ rounds/second is the number of revolutions per second of the propeller.   

Moreover, the $\kappa-\omega$ SST model is coupled with a $\gamma-Re_{\theta_t}$ transition model, which simulates the laminar-turbulent transition in our simulation setting and was proposed in~\cite{langtry2009correlation}. 
The initial and boundary conditions for the additional variables $\kappa$ and $\omega$, $\gamma$ and $Re_{\theta_t}$ are set in the \fom{} following standard formulas and rules.
In particular, the turbulent kinetic energy $\kappa$ and the turbulence specific dissipation rate $\omega$ are initialized to:
$$\kappa = \dfrac{3}{2}(I \,u_0)^2\,, \omega = \frac{\kappa^{0.5}}{C_{\mu}^{0.5}L}$$
in all the domain, where the $I$ is the turbulence intensity, that is $\sim 9 \%$ in our case, $C_{\mu}=0.09$ and $L=D$ is the reference length scale.
The intermittency variable $\gamma$ and the transition momentum thickness Reynolds number $Re_{\theta_t}$ are initialized to:
$$\gamma=1\,, Re_{\theta_t}= \begin{cases}
    1175.51-589.428I+\frac{0.2196}{I^2} \text{ if }I<1.3\, ,\\
    \frac{331.5}{(I-0.5658)^{0.671}} \text{ if }I>1.3\,.
    \end{cases}$$
For what concerns the boundary conditions, we have for all variables a fixed-value condition at the inlet, and a zero-gradient condition at the outlet and at the walls.

\RA{As stated before, in this paper we adopt} the finite volume discretization technique\RA{, which }consists in performing a polyhedral discretization of the computational domain, where each finite volume is called \emph{control volume}. After that, the equations in \eqref{U-RANS-turb} are written in integrated on each control volume of the domain interest. The divergence theorem is then used to convert the volume integrals
to surface integrals, which are finally discretized as sums of the fluxes at the boundary faces of each control volume.

\subsection{\mrf{} approach}
\label{subsec:mrf}
For what concerns the motion of the propeller, the rotation of the propeller in \emph{OpenFOAM} is obtained using the \emph{Moving Reference Frame} (\mrf{}) approach. The \mrf{} technique is a steady-state method widely spread in industrial \cfd{} problems that involve rotating parts. The principle of this technique is the creation of a thin volumetric region of mesh cells around the rotating body during the meshing phase. This region, namely the \emph{MRF zone}, is in our test case the smallest cylinder which surrounds the propeller's blades ($\Gamma_1$ in Figure \ref{fig:domain}). During the simulation, the MRF zone is rotated about the axis of the body and the body is kept stationary. The simulation is performed until a steady state is reached, such that the thrust and the torque forces reach a stationary value. 

This stage aims to build a \cfd{} model which provides accurate results compared to the experimental open-water tests.
In particular, our metrics for accuracy measurement are the relative error on the values of the thrust and torque coefficients ($k_T$ and $k_Q$). These adimensional numbers are characteristic features in the fluid dynamics of marine propellers and are defined as follows:
\begin{equation}
k_T = \dfrac{T}{\rho n^2 D^4}, \quad k_Q = \dfrac{Q}{\rho n^2 D^5}.
    \label{eq:coeffs_def}
\end{equation}
In expression \eqref{eq:coeffs_def}, $\rho$ is the fluid density.
$T$ and $Q$ are the thrust and the torque force experimented on the blades of our propeller. These forces are evaluated as the sum of a pressure-based and a viscous contribution, as follows:
    \begin{equation}
        \begin{split}
            & T = T_{\text{pressure}}+T_{\text{viscous}} = \rho \left(\int_{\Gamma_{\text{blades}}} p \mathbf{n} dA + \int_{\Gamma_{\text{blades}}} \mathbf{\Sigma n} dA \right);\\
            & Q = Q_{\text{pressure}}+Q_{\text{viscous}} = \rho \left(\int_{\Gamma_{\text{blades}}} p \mathbf{n} \times \mathbf{r} dA + \int_{\Gamma_{\text{blades}}} \mathbf{\Sigma n} \times \mathbf{r} dA \right).
        \end{split}
        \label{eq:T_and_Q}
    \end{equation}
    In expression \eqref{eq:T_and_Q}, $\Gamma_{\text{blades}}$ indicated the blades' surface, $\mathbf{r}=(x,y,z)$ is the position vector, $\mathbf{n}$ is the unitary normal vector to the surface, $\mathbf{\Sigma}$ is the wall-shear stress tensor.

In the project presented in this paper, the \fom{} performed in \emph{OpenFOAM} reached an accuracy of $\sim 1 \%$  and of $\sim 3 \%$ for what concerns $k_T$ and $k_Q$, respectively.
Moreover, the time needed to perform one high-fidelity simulation is $24-48$ hours in a parallel setting \footnote{The \fom{} simulations are performed using 55 processor cores on SISSA HPC cluster Ulysses (200 TFLOPS, 2TB RAM, 7000 cores).} The variance in the simulation time is due to the different number of iterations needed to reach a stable regime in all the computed simulations.

\section{Non-intrusive Reduced Order Model}
\label{sec:roms}
This Section is dedicated to the explanation of the theory behind the model order reduction technique here employed to reduce the computational effort. This technique belongs to the framework of \emph{non-intrusive reduced order methods}, i.e. approaches exploiting the information provided by the high-fidelity simulations. In these kinds of models, the \cfd{} governing equations are only used at the full order level to perform the offline simulations, but not at the reduced order level.
Non-intrusive reduced order models are composed of two fundamental stages, the \emph{reduction} and the \emph{approximation} step, described in the following paragraphs. The techniques here explained are implemented in the Python package \verb+EZyRB+ \cite{DemoTezzeleRozza2018EZyRB, ezyrb}.

\paragraph{Reduction}
The first step consists of the compression of the original matrix of high-fidelity solutions, into a matrix of reduced dimension. In our case, we consider as snapshots the pressure and wall shear stress at the final time instant of each offline simulation, i.e. when a steady state has been reached. The fields are evaluated not on all the points of the computational mesh, but only on the points of the blades. 

It is important to remark that the flow fields are evaluated only on the blades because the propeller efficiency only depends on the fluid-dynamics behaviour at the blades. In fact, the efficiency is evaluated as:
\begin{equation}
    \eta_{\text{propeller}}=\dfrac{T}{Q} \dfrac{u_0}{2 \pi n},
    \label{eq:eta}
\end{equation}
where $T$ and $Q$ only depend on the fields on the blades, as can be evinced from \eqref{eq:T_and_Q}.

Thus, each $i-$th snapshot is a flow field evaluated on the deformed blades corresponding to the $i-$th set of deformation parameters taken into account. For instance, if we consider the generic field $s$, we have the following matrix of snapshots:
$$\mathbf{S} = \begin{bmatrix}
    \vert & \vert & &\vert\\
    \mathbf{s}_1(\mathbf{x}) & \mathbf{s}_2(\mathbf{x}) & \dots & \mathbf{s}_\nsnaps{}(\mathbf{x})\\
    \vert & \vert & &\vert
\end{bmatrix}\in \mathbb{R}^{N_{\text{dof}} \times \nsnaps{}},$$
where $N_{\text{dof}}$ is:
\begin{itemize}
    \item exactly the number of cells of the blades in the computational mesh, that is $\sim \num{28e4}$ in the first type of \rom{} we consider (the \emph{\textbf{standard} ROM}). In this case, the snapshots are evaluated at the cells' centers;
    \item a certain number of quadrature points retained on the blades, in our case $\num{12e4}$, since we consider a grid of $100\times 100$ nodes for the back and the face of each blade, in the second type of \rom{} (the \emph{\textbf{fast} \rom{}}). In this case, the snapshots are evaluated on the quadrature points. 
\end{itemize}
The basic principles of the two theories will be explored in Subsections \ref{rom_1} and \ref{rom_2}.

The reduction technique adopted in this work is the Proper Orthogonal Decomposition (\podec{}). This technique is based on the projection of the snapshots into a reduced space, spanned by a limited number of the so-called \emph{modes}, which are computed directly starting from the snapshots matrix in the offline stage.
Each \emph{reduced} snapshot $\mathbf{s}_i$ can be approximated as a linear combination of the modes: $$\mathbf{s}_i \simeq \sum_{j=1}^L a_j \boldsymbol{\phi}_j,$$
where $\{\mathbf{\phi}_j\}_{j=1}^L$, are the modes, $L \ll N_{\text{dof}}$ is the reduced dimension, that has to be established a priori. $\{a_j\}_{j=1}^L$ are the reduced coefficients associated to the modes.
The \podec{} modes can be evaluated using a Singular Value Decomposition technique (SVD) or via the correlation matrix.
In the first case, for instance, the snapshots matrix is decomposed in $\mathbf{S}=\mathbf{U \Sigma} \mathbf{V}^T$, where the columns of $\mathbf{U} \in \mathbb{R}^{N_{\text{dof}} \times L}$ are the \podec{} modes. Thus, the reduced coefficients can be computed as $\mathbf{U}^T \mathbf{S}$.

\paragraph{Approximation}
In the second step, the goal is to predict the reduced coefficients associated with unknown values of the parameters, which do not belong to the original dataset.
In our case, the snapshots are $\mathbf{s}_i= \mathbf{s}(\boldsymbol{\mu}_i)$, $i=1,\dots,M$ and the parameters are $\boldsymbol{\mu}_i \in \mathbb{R}^p$, $p=4$.
The goal here is to evaluate $\mathbf{s}(\boldsymbol{\mu}^*)$, i.e. the pressure field for deformation parameters that are not in our initial set of parameters.

Different techniques can be employed to reach this task. Here we consider one of the most used techniques, the \rbf{} interpolation. The \rbf{} allows representing our field at unknown parameter $\boldsymbol{\mu}^*$. In particular, the unknown coefficients $\mathbf{a}^*(\boldsymbol{\mu})$ can be computed in the following way:
$$\mathbf{a}(\boldsymbol{\mu}^*)=\sum_{i=1}^M \omega_i \phi(\| \boldsymbol{\mu^*} - \boldsymbol{\mu}_i \|), $$
where $\phi(\| \boldsymbol{\mu^*} - \boldsymbol{\mu}_i \|)$ are the radial basis functions with center $\boldsymbol{\mu}_i$ and weight $\omega_i$. The weights are found from the conditions:

$$\mathbf{a}(\boldsymbol{\mu}_j) = \sum_{i=1}^M \omega_i \phi(\| \boldsymbol{\mu}_j -\boldsymbol{\mu}_i\|)\,, j=1, \dots, M\,.$$

In this work, radial basis functions of multiquadric shape, of expression $\phi(r)=\sqrt{1+(\varepsilon r)^2}$, with $r=\|\boldsymbol{\mu}^* - \boldsymbol{\mu}_i\|$.
Alternative approaches that can be used as \emph{approximation part} of \rom{}s are the Gaussian Process Regression (GPR) \cite{williams2006gaussian}, or the K-Neighbors Regression (KNR).

\subsection{Standard \rom{}: full mesh}
\label{rom_1}

In the standard-\rom{} model we introduce, we consider as snapshots the pressure and wall shear stress fields evaluated on all the mesh points of the blades. 

The question that automatically arises is: \emph{What is the procedure to predict the efficiency of a propeller starting from the deformation parameters $\boldsymbol{\mu}^*$?}

The efficiency of a propeller can be predicted by performing the following steps:
\begin{enumerate}
    \item \label{def-blade-std} deform the original blades according to the deformation parameters $\boldsymbol{\mu}^*$;
    \item \label{def-mesh-std} deform the blades' points using a \rbf{} technique and triangulate the mesh;
    \item \label{na-std} compute the normal unitary vectors to the cells and their areas, that are necessary for the computation of thrust and torque forces;
    \item \label{rom-std} exploit the \emph{standard} \rom{} to predict the pressure and wall-shear stress fields on the blades;
    \item \label{eff-std} compute the thrust and torque forces approximating the integrals in \eqref{eq:T_and_Q}, and then the efficiency (following expression \eqref{eq:eta}).
\end{enumerate}

\textbf{Remark} \\ We specify here that in step \ref{def-mesh-std} a triangulation is applied on the deformed mesh to make straightforward the computation of the normals and areas to the cells. \\

The integrals \ref{eq:T_and_Q} are approximated as discrete sums on all the mesh cells, as follows:
\begin{equation}
\begin{split}   
    &T_{\text{approx}} = \rho \sum_{c=1}^{N_{\text{dof}}}\left(p_c \mathbf{n}_{c} a_{c} + \boldsymbol{\Sigma}_c \mathbf{n}_{c} a_{c} \right)\,;\\
    &Q_{\text{approx}} = \rho \sum_{c=1}^{N_{\text{dof}}} \left( 
    p_c \mathbf{n}_{c} \times \mathbf{x}_{c} a_{c} + \boldsymbol{\Sigma}_{c} \mathbf{n}_{c} \times \mathbf{x}_{c}  a_{c}\right)\,,
    \end{split}
\end{equation}
where $\mathbf{x}_c$ is the position of the centre of cell $c$, the approximated reduced fields are $p_c=p_{\text{rom}}(\mathbf{x_c})$, $\boldsymbol{\Sigma}_c=\boldsymbol{\Sigma}_{\text{rom}}(\mathbf{x_c})$, $\mathbf{n}_c$ and $a_c$ are the normal to cell $c$ and its area.

The critical issue of this analysis is the computational time required to go through steps \ref{def-blade-std}-\ref{eff-std}. In particular, the most computationally expensive step is \ref{def-mesh-std} and the total time to predict the efficiency of a deformed propeller is $\sim 6$ minutes.

\subsection{Fast \rom{}: quadrature points}
\label{rom_2}
We here consider the second type of \rom{}, named \emph{fast} \rom{}.
In this case, the snapshots are evaluated on the coordinates of a pre-defined number $N_{\text{quadrature}}$ of Gauss quadrature nodes on the blades.

\RA{The 3-dimensional quadrature points are found mapping the 2D cartesian Gauss-Legendre quadrature nodes into the UV local reference system on the NURBS surfaces, and then mapping the UV coordinates in a 3D cartesian domain. In particular, we considered as different NURBS the suction and pressure side of all blades.}

\RA{\textbf{Remark} \\ 
The number of quadrature nodes on the blades has been selected after conducting a sensitivity study. This study is conducted on the predictions of the thrust forces and torque momentum. In particular, we measured the accuracy of the prediction made considering quadrature formulas with respect to the FOM value, for different degrees of quadrature nodes.
The results are represented in Figure \ref{fig:sensitivity-quadrature}. We can notice that the accuracy converges on steady values for a degree $\geq 100$. For this reason, we selected $100$ quadrature nodes on each face of the blades.

\begin{figure}[h!]
    \centering
    \includegraphics[width=\textwidth]{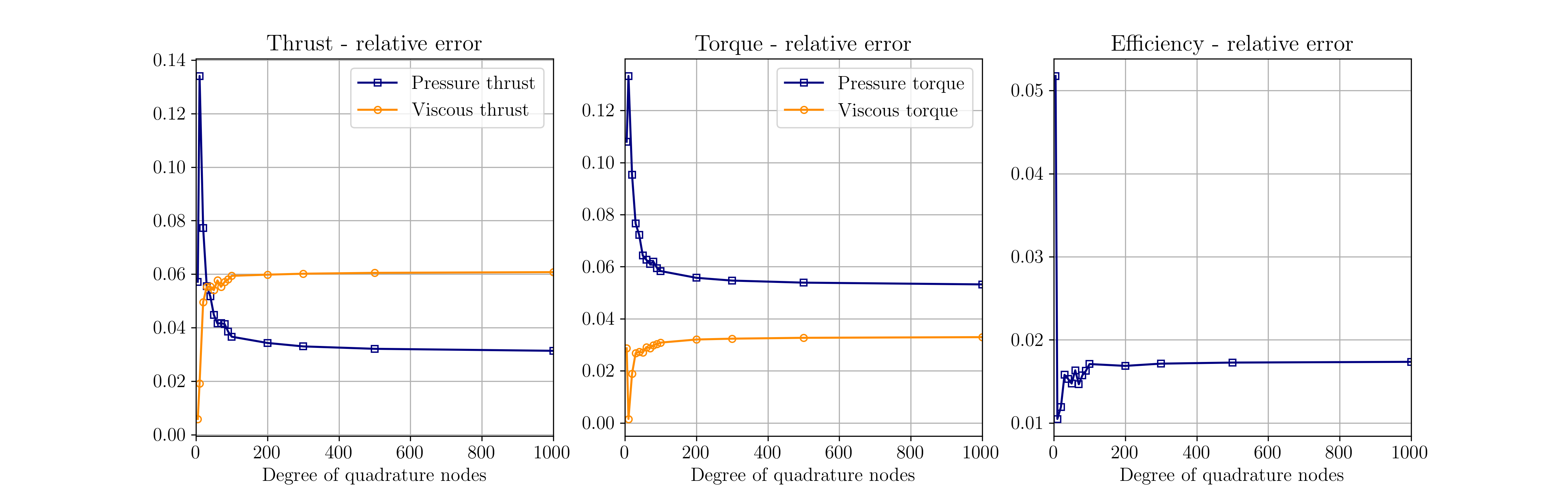}
    \caption{Accuracy of predicting forces, momentum and efficiency using quadrature formulas, for different number of quadrature samples. The accuracy is measured as a relative error with respect to the reference FOM value.}
    \label{fig:sensitivity-quadrature}
\end{figure}

}

\begin{figure}
    \centering
    \subfloat[Undeformed blades' quadrature nodes]{\includegraphics[width=0.3\textwidth]{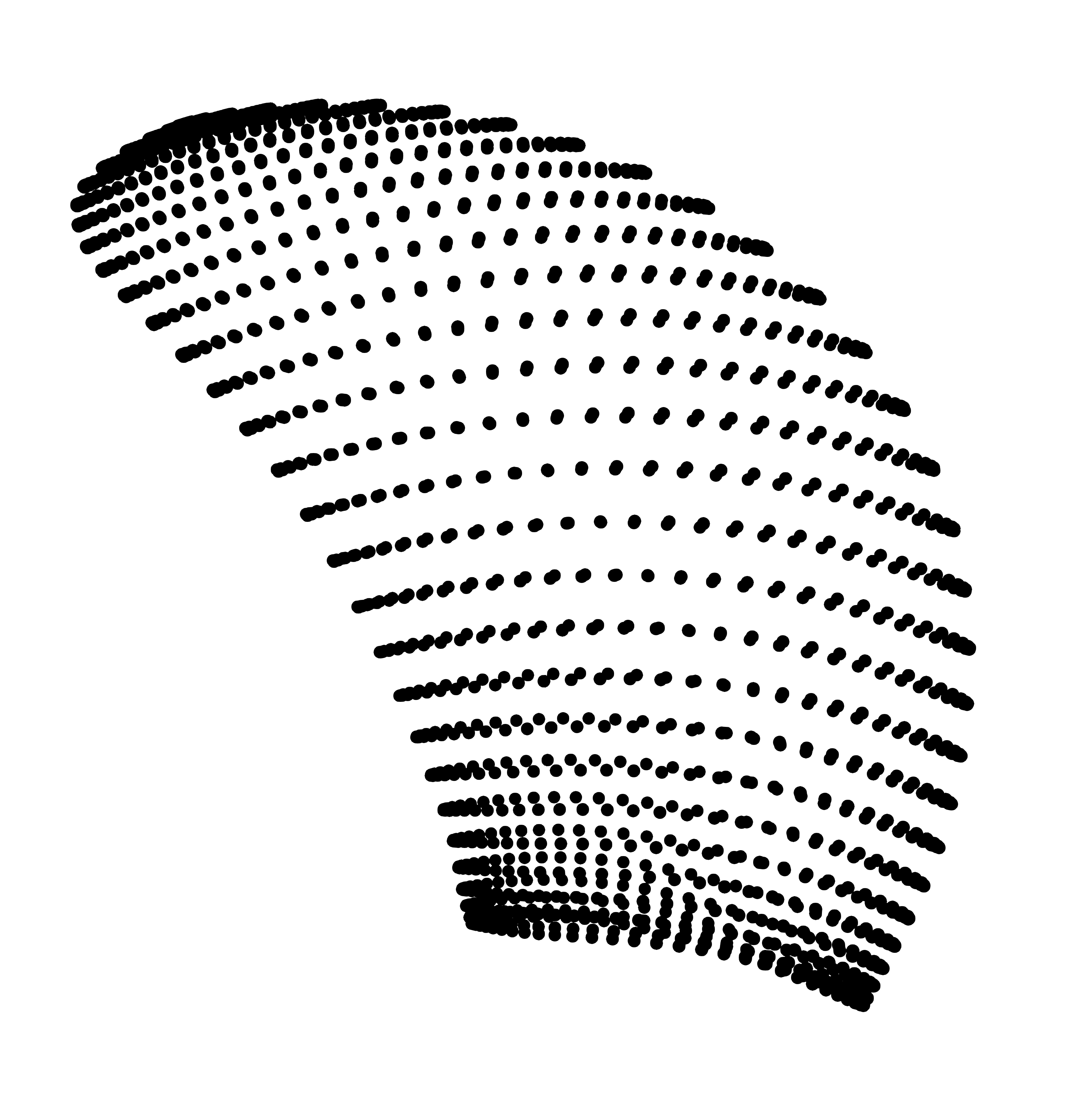}}
    \subfloat[Deformed blades' quadrature nodes]{\includegraphics[width=0.3\textwidth]{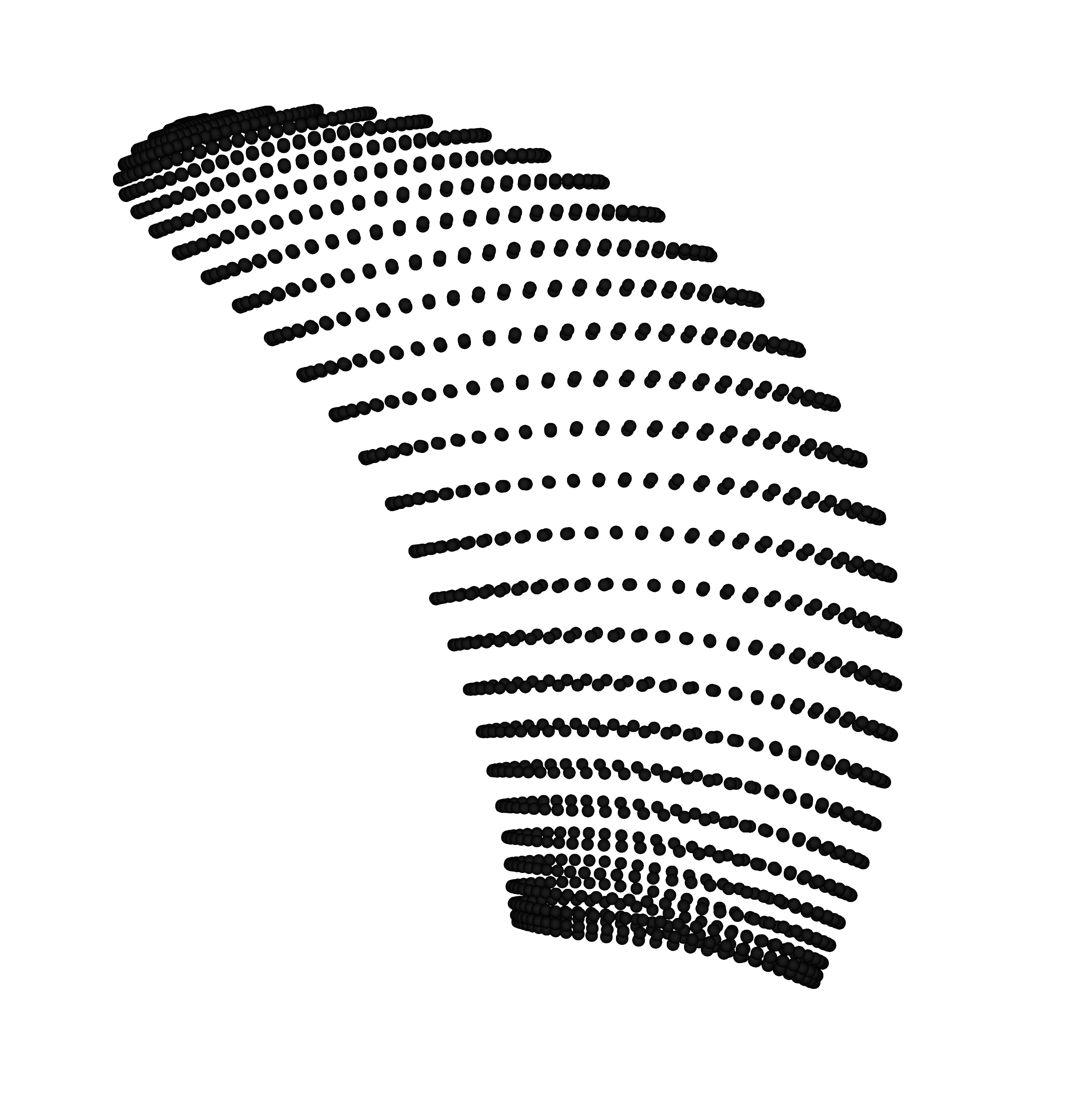}}
    \caption{Quadrature nodes on undeformed and deformed blades. In this picture, $N_{\text{quadrature}}=10800$.}
    \label{fig:quadrature-nodes}
\end{figure}

In Figure \ref{fig:quadrature-nodes}, two representations of the Gauss quadrature nodes on the undeformed and deformed blades are displayed. 

Therefore, there are in the offline stage two additional steps with respect to the standard \rom{}:
\begin{itemize}
    \item the computation of the quadrature nodes on all the deformed blades in the dataset;
    \item the evaluation of the snapshots on the quadrature nodes through \emph{interpolation} or \emph{regression} techniques. Here we consider a K-Neighbor Regression technique based on 5 nearest neighbors, where the fields on the quadrature nodes are predicted by local interpolation of the targets associated with the nearest neighbors in the training set.
    \end{itemize}
Moreover, a \rom{} is also exploited not only for the prediction of the pressure and wall shear stress fields but also for the computation of the normal vectors to the blades' surfaces.

We illustrate here the complete procedure to predict the efficiency of a deformed propeller with parameters $\boldsymbol{\mu}^*$:
\begin{enumerate}
    \item \label{def-blade-fast} deform the initial blades according to parameters, as in \ref{def-blade-std};
    \item \label{quad-fast} generate a number $N_{\text{quadrature}}$ of Gauss quadrature nodes on the deformed blades;
    \item \label{normal-fast} predict the normal vectors exploiting \rom{} and compute the jacobians of the deformation;
    \item \label{rom-fast} exploit the \rom{}s to predict the pressure and wall shear stress on the blades' quadrature points;
    \item \label{eff-fast} evaluate the torque and thrust forces, and consequently the efficiency.    
\end{enumerate}

The technique here presented, which --- to the best of our knowledge --- is novel for this application, allows for a reduction in the computational time.
Indeed, the efficiency evaluation for a single propeller takes approximately $8-15$ seconds \footnote{Both the \rom{}s are performed on an Intel(R) Core(TM) i5-4570 CPU @ 3.20GHz 16GB RAM on only one processor core.}. The reason for this gain is that the mesh generation, i.e. the expensive step \ref{def-mesh-std} in standard \rom{}, is here replaced by step \ref{quad-fast}, which is much faster than the mesh generation.

The most important difference between the two approaches is the computation of the forces: in the standard case, we approximate the integrals in \eqref{eq:T_and_Q} with discrete sums on all the mesh cells; in the second case, the quadrature formulas on the Gauss nodes are considered to approximate the integrals.
However, the use of integration formulas to compute the forces does not lead to a minor accuracy in the reconstruction of the efficiency.

\subsection{Results and comparison of \rom{}s}
\label{subsec:results-roms}
As specified in the previous paragraphs, the reduced order model can be built considering different techniques for the reduction and approximation stages.
In this work, we always exploit the \podec{} as a reduction method, since it is a consolidated technique in the field of industrial applications. On the other side, the approximation technique has been chosen as the result of an analysis of the accuracy of the resulting \rom{}. In particular, the \rbf{}, \gpr{} and \knr{} approximating techniques are here compared.
As for the \rbf{}, the thin plate spline is chosen as basis functions; in the \knr{} the regression is implemented considering $k=5$ and uniform weights.
The accuracy is here measured in terms of \emph{k-fold cross validation} error. It is computed by splitting the database into $k$ consecutive folds (we chose $k=10$ in our case), and each fold is used once as validation while the $k - 1$ remaining folds form the training set. We remind that the database is composed of $\nsnaps{}=216$ snapshots, corresponding to the initially deformed blades (details about the sampling are provided in~\ref{subsec:param}.

The mean values of the \emph{k-fold} errors for the standard and fast reduced order models are reported in Tables \ref{tab:rom_analysis_std} and \ref{tab:rom_analysis_fast}.
The results show that the \gpr{} is the best technique for the standard approach, whereas the \rbf{} provides the best results for the fast \rom{}. The best methods are chosen for the efficiency prediction in the shape optimization in order to obtain good results in the \fom{}-\rom{} validation.

\begin{figure}[htbp]
    \centering
    \subfloat[\fom{}]{\includegraphics[width=0.28\textwidth, trim={30cm, 0cm, 30cm, 0cm}, clip]{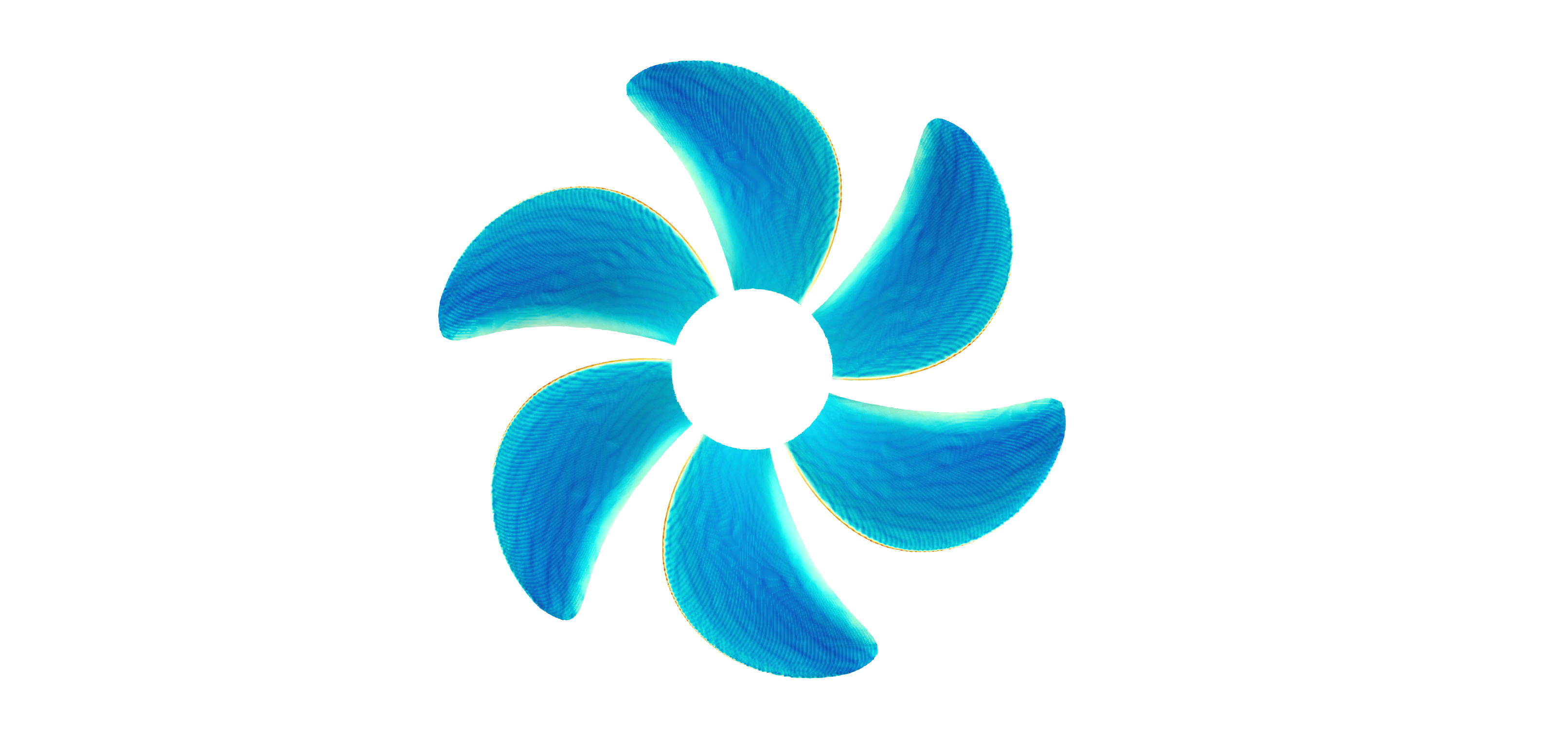}\label{fom-p}}
    \subfloat[Standard \rom{}]{\includegraphics[width=0.28\textwidth, trim={30cm, 0cm, 30cm, 0cm}, clip]{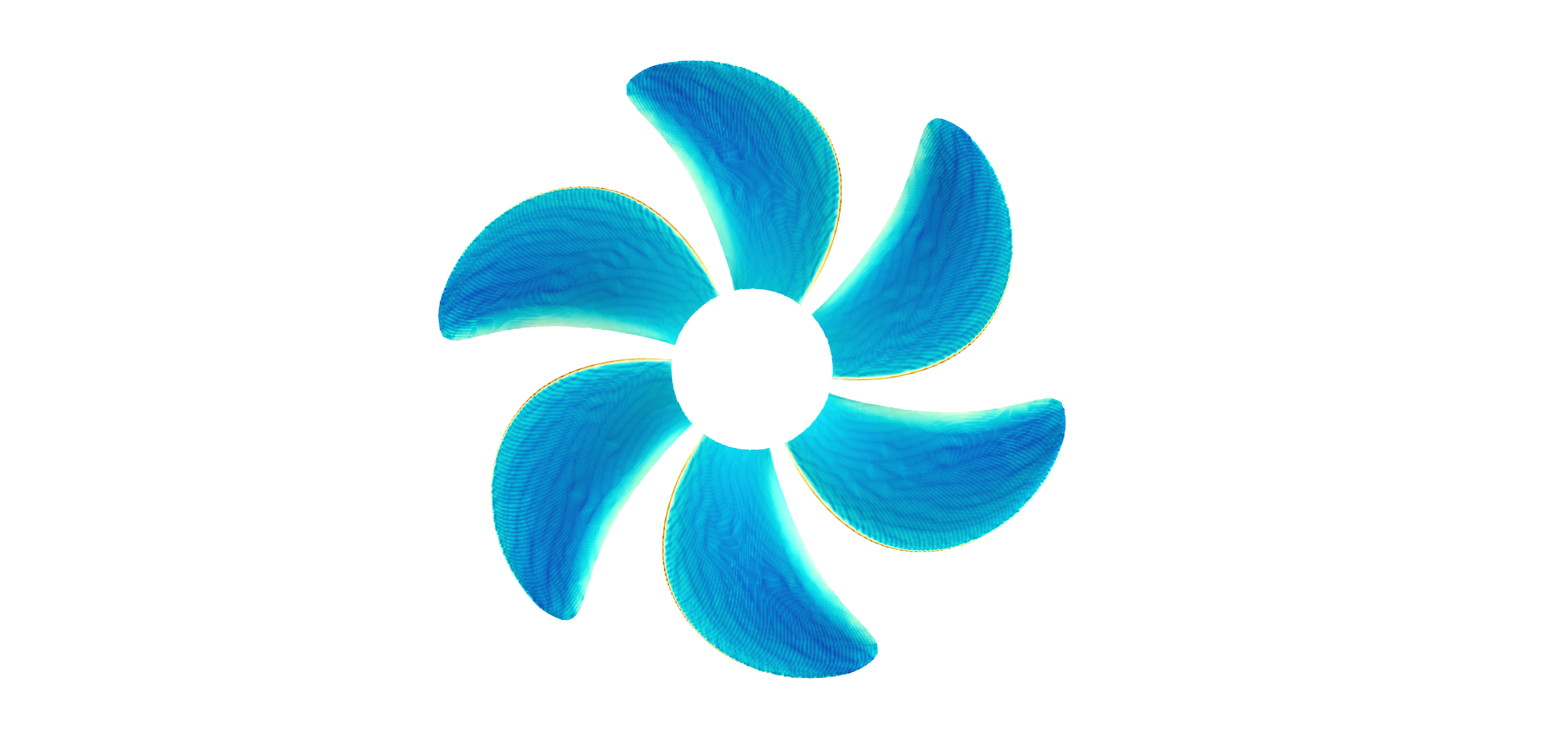}\label{rom-std-p}}
    \subfloat[Fast \rom{}]{\includegraphics[width=0.41\textwidth, trim={30cm, 0cm, 6cm, 0cm}, clip]{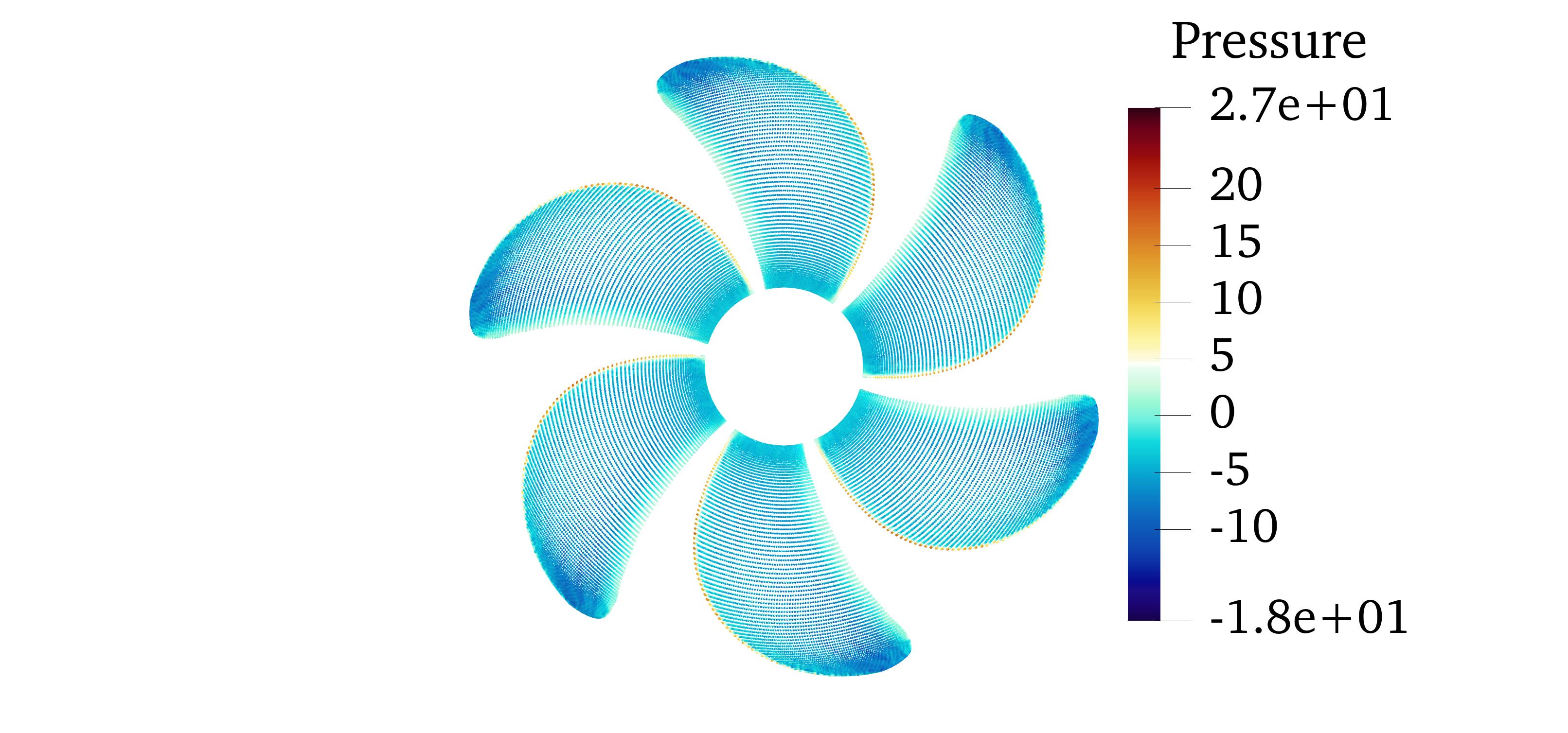}\label{rom-fast-p}}\\
    \subfloat[\fom{}]{\includegraphics[width=0.28\textwidth, trim={30cm, 0cm, 30cm, 0cm}, clip]{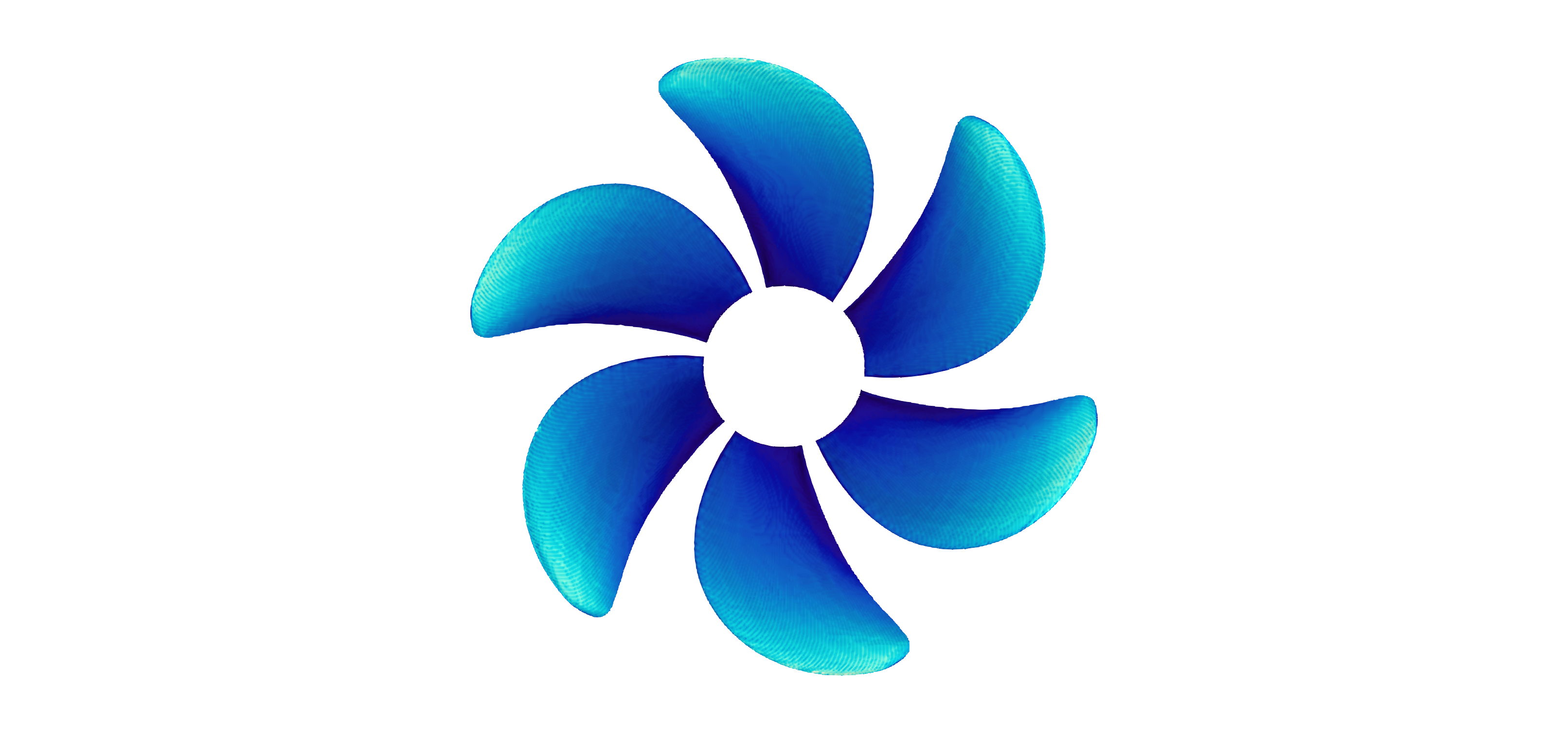}\label{fom-shear}}
    \subfloat[Standard \rom{}]{\includegraphics[width=0.28\textwidth, trim={30cm, 0cm, 30cm, 0cm}, clip]{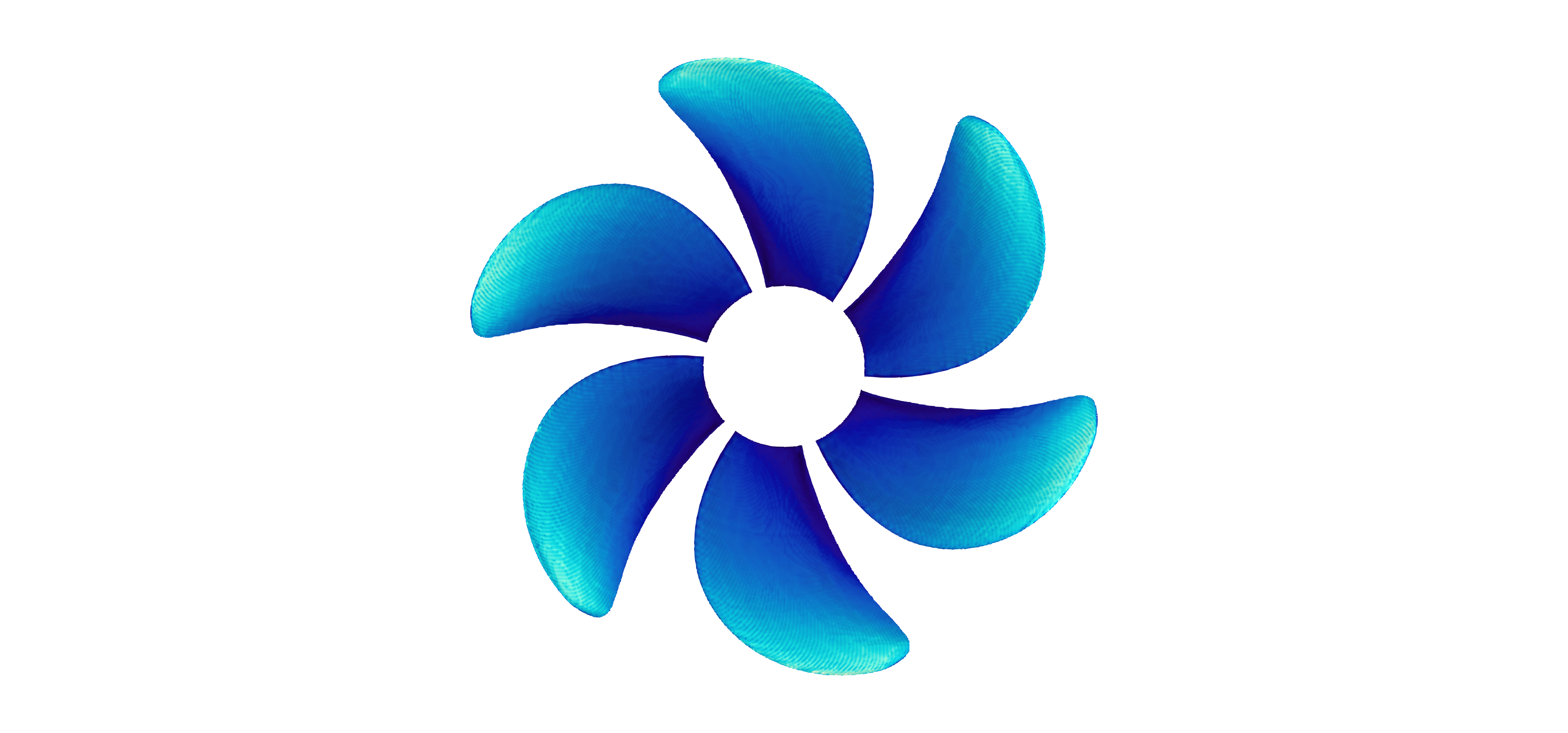}\label{rom-std-shear}}
    \subfloat[Fast \rom{}]{\includegraphics[width=0.41\textwidth, trim={30cm, 0cm, 5.7cm, 0cm}, clip]{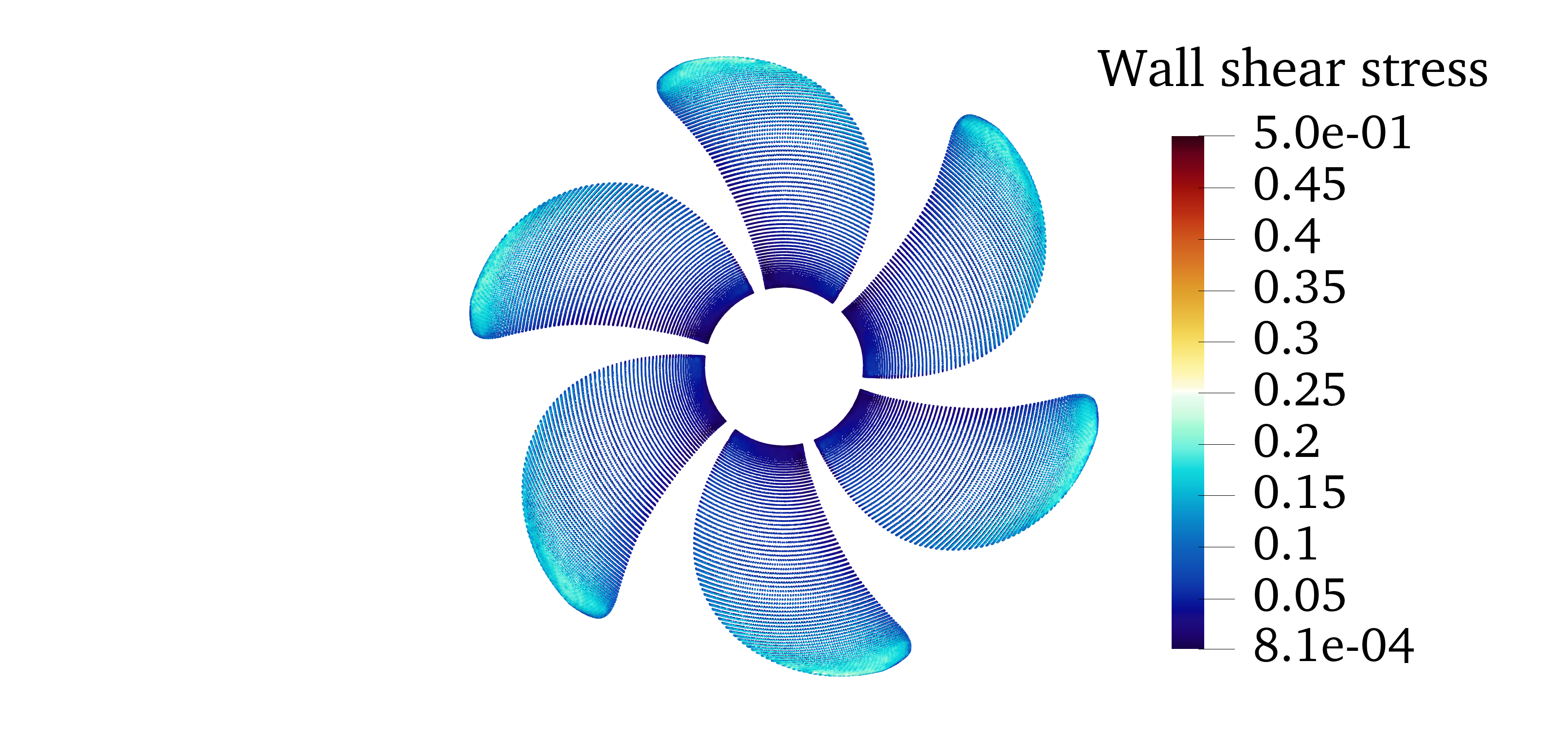}\label{rom-fast-shear}}
    \caption{High-fidelity and reduced order fields corresponding to the initial blade: relative pressure in the first row and wall shear stress magnitude in the second row.}
    \label{fig:fom-rom}
\end{figure}

Figure \ref{fig:fom-rom} displays a comparison between the full-order fields acting on the blade and the reconstructed reduced-order fields. In the case of the fast \rom{} approach (Figure \ref{rom-fast-p} and \ref{rom-fast-shear}), the fields are computed and represented on the Gauss quadrature nodes.
The Figure shows a good agreement between the reconstructed fields and their original high-fidelity counterpart.

\section{Optimization process}
\label{sec:opt}
This Section is dedicated to the shape optimization process.  
In particular, Subsection \ref{subsec:genetic} briefly recalls the theory on the genetic algorithm, here exploited for optimization, and Subsection \ref{subsec:gradient} is dedicated to a comparison among the results obtained from the genetic optimization and those obtained using a gradient-based approach.

\subsection{Genetic optimization}
\label{subsec:genetic}

This Subsection is dedicated to the explanation of the Genetic Algorithm (GA) used for the optimization step. It was introduced for the first time in \cite{holland1973genetic} and it is inspired by the evolution mechanism.

\begin{figure}
    \centering
    \includegraphics[width=0.8\textwidth]{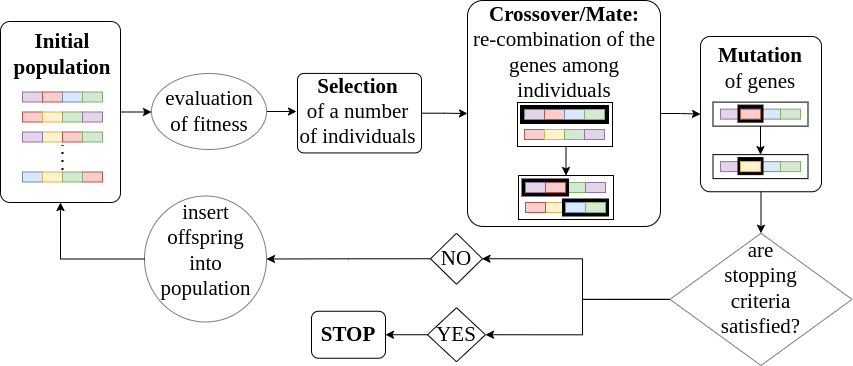}
    \caption{Schematic representation of the genetic algorithm.}
    \label{fig:genetic}
\end{figure}

The general scheme of how it works is displayed in Figure \ref{fig:genetic}. In particular, it starts from an initial population of $N_{\text{pop}}$ individuals and the \emph{fitness}, i.e. the objective function we want to minimize/maximize, is evaluated for all the individuals. Then, the algorithm iteratively performs three steps: \emph{selection} of the best individuals, \emph{crossover}, recombination of the genes of the individuals, and \emph{mutation} of the genes, allowing the evaluation of new individuals.

The GA is computed in Python making use of the \emph{Deap} evolutionary computational framework~\cite{de2012deap, fortin2012deap} and the following hyper-parameters are considered:
\begin{itemize}
    \item the initial population is composed by $30$ and $150$ individuals for the standard and fast \rom{}s, respectively;
    \item in the \emph{mate} stage, a one-point crossover is performed, that modifies in-place the input individuals;
    \item in the \emph{mutation} step, a Gaussian mutation of mean $\mu=1$ and standard deviation $\sigma=0.1$ is applied to the input individual. The independent probability for each attribute to be mutated is set to $0.5$ and $0.8$, when the standard and fast \rom{}s are exploited, respectively;
    \item the \emph{evolutionary algorithm} employed in the genetic process is the $(\mu + \lambda)$ algorithm, which selects the $\mu$ best individuals for the next generation ($\mu=5$ for the standard \rom{}, $50$ in the fast approach) and produces $\lambda$ children at each generation ($\lambda=10$ in the standard approach, $80$ in the fast one). The probabilities that an offspring is produced by crossover and by mutation are $0.4$ and $0.5$, respectively.
    \item the previous steps are repeated for $10$ and $20$ generations for the standard and fast cases, respectively.
\end{itemize}

\textbf{Remark} \\ In this paper, we chose the genetic algorithm for optimization instead of a gradient-based algorithm since it allows us to evaluate a larger number of individuals, in our case of deformed shapes, and avoids getting stuck into local minima.
\\

In our particular test case, the individuals are the parameters, so each individual $\mathbf{i}$ has 4 genes, corresponding to the deformation rates of the geometrical parameters.
The fitness function $f(\mathbf{i})$ we want to maximize is the efficiency. We consider two different optimization processes:
\begin{itemize}
    \item basic optimization, where the efficiency $\eta_{\text{propeller}}$ is the fitness;
    \item constrained optimization, where a series of geometrical and physical constraints applies to the optimization process. In this case, we consider a penalized fitness function: 
    \begin{equation}
    f(\mathbf{i}) = \eta_{\text{propeller}} (\mathbf{i}) - \sum_{j=1}^{N_{\text{constraints}}} w_j \, \text{penalty}(\mathbf{i})_j,
    \label{fitness_penalized}
    \end{equation}
    where $\text{penalty}(\mathbf{i})_j$ and $w_j$ are the penalty and the weight associated with the $j-$th constraint, $N_{\text{constraints}}$ is the number of constraints taken into account.
\end{itemize}
We specify here only the generic formulation to apply constraints in the optimization procedure since
the constraints imposed in this work are protected by Non-Disclosure Agreements. It derives that a punctual discussion about the geometrical features of the optimal blade can not be pursued; we consider the constrained optimization as an additional test \RB{with} a different objective function, limiting our considerations to the precision comparison between \rom{} and \fom{}.

\paragraph{Results of genetic optimization}
Table \ref{tab:genetic} displays the best individual obtained from the genetic algorithm (GA) in the standard and fast \rom{} approaches.

In particular, we measure the results of the optimization processes in terms of absolute increases of the percentage efficiency with respect to the starting propeller, i.e. the one corresponding to parameter $\mu = \begin{bmatrix}    1&1&1&1   \end{bmatrix}$. Table \ref{tab:genetic} reports the following metrics:
\[
\Delta_{\text{\fom{}}} = \eta^{\text{opt}}_{\text{\fom{}}} - \eta^{\text{start}}_{\text{\fom{}}}, \quad
\Delta_{\text{\rom{}}} =  \eta^{\text{opt}}_{\text{\rom{}}} - \eta^{\text{start}}_{\text{\rom{}}},
\]
where $\eta^{\text{start}}_{\text{\fom{}}}$ and $\eta^{\text{start}}_{\text{\rom{}}}$ are the \fom{} and \rom{} efficiencies of the starting propeller, whereas $\eta^{\text{opt}}_{\text{\fom{}}}$ and $\eta^{\text{opt}}_{\text{\rom{}}}$ are the efficiencies for the optimal propeller.
Moreover, we underline here that the efficiency of the starting propeller is captured with an absolute error of $0.01 \%$ in the standard \rom{} approach and of $0.33 \%$ in the fast \rom{} approach.

Table \ref{tab:genetic} shows that the standard and fast optimization processes lead to similar values of the final optimal efficiency, although the algorithms converge to different results in terms of individuals in the parameter space.

In general, the optimized propeller obtained with an unconstrained algorithm is characterized by higher efficiency with respect to the one resulting from the constrained process. Indeed, as can be seen from the Table, the unconstrained optimization processes lead to blades with small thickness, since it would result in lower torque values, and hence to more efficient propellers.
On the other hand, thin blades can reduce the robustness of the propeller's structure and, hence, make the propeller more vulnerable to external stress. This is the reason why geometrical constraints (on the chord length and on the thickness, for example) are required by industries to ensure the robustness of the propeller's mechanical structure.

For what concerns the accuracy in the prediction of the efficiency, Table \ref{tab:genetic}, together with Tables \ref{tab:rom_analysis_std} and \ref{tab:rom_analysis_fast}, show that the standard technique is more accurate, with an error behind the $0.1 \%$ on the efficiency value.
However, the fast algorithm allows for reaching parametric points with comparable full-order efficiencies in a significantly reduced amount of time.
Indeed, as already pointed out in Sections \ref{rom-std} \ref{rom-fast}, the evaluation of the efficiency for a single individual lasts about $6$ minutes exploiting the standard \rom{}, while it lasts less than $15$ seconds exploiting the fast approach.

\begin{figure}
    \centering
    \subfloat[]{\includegraphics[width=0.35\textwidth]{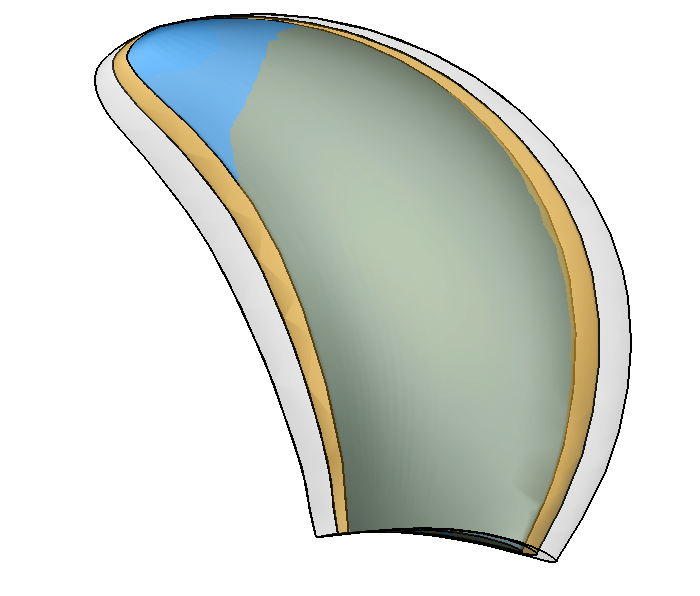}\label{subfig:stand-opt}}
    \subfloat[]{\includegraphics[width=0.35\textwidth]{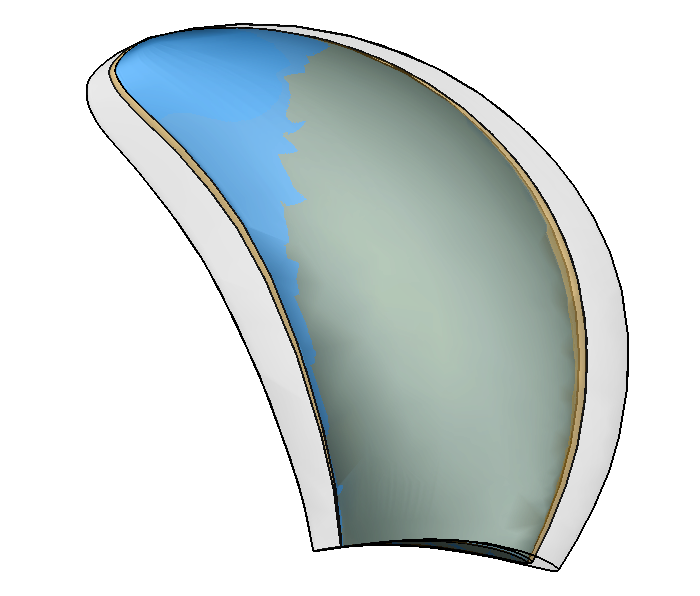}\label{subfig:fast-opt}}
    \caption{Graphical representation of the final optimal shapes obtained from the standard \protect \subref{subfig:stand-opt} and the fast \protect \subref{subfig:fast-opt} optimization processes. The Figure displays the starting blade in gray, the optimal blade from the unconstrained optimization process in light blue, and the optimal blade from the constrained optimization in orange.}
    \label{fig:optimals}
\end{figure}

Another consideration that can be withdrawn is that the two techniques do not converge to the same optimal parametric points.
The reason for this fact is that the objective function is computed following a different approximated procedure. The same holds true for the penalty terms related to forces' constraints in the constrained optimization. The optimal blades are graphically represented in Figure \ref{fig:optimals} and compared in shape with the original starting blade of the pipeline (colored in gray).

\subsection{A comparison with gradient-based algorithms}
\label{subsec:gradient}

In this paper, the authors chose to adopt a genetic algorithm for the optimization process, since it allows the evaluation of a large number of individuals, increasing the probability that the global optimum is reached. On the other hand, gradient-based algorithms suffer from high sensitivity to the initial guess of the optimization and this can lead the algorithm to freeze into the closest local optimum.
Moreover, the genetic process is characterized by a large number hyper-parameters, specified in Subsection \ref{subsec:genetic}, and changing those parameters allowed us to do a large number of optimization experiments and to choose the values resulting in the highest objective function.

Here we show different experiments done using a gradient-based algorithm varying the initial guess of the optimization process. This part focuses on the results obtained using two different techniques: the \emph{conjugate gradient} (\emph{CG}) method, first described in~\cite{hestenes1952methods}, and the \emph{L-BFGS-B} method, namely the \emph{Limited-Memory} version of the bound-constrained \emph{BFGS} (\emph{Broyden–Fletcher–Goldfarb–Shanno} algorithm, treated in~\cite{broyden1970convergence, robitaille1970quasi, byrd1995limited}.

In this Subsection, we focus on the \textbf{constrained fast optimization} and made a comparison between the results obtained with the above-cited gradient-based methods and the genetic algorithm.
The results are reported in Table \ref{tab:grad} for different starting points and for the CG and L-BFGS-B approaches.
We specify that the fitness increments $\Delta_{\text{ROM}}$ are always measured with respect to the initial propeller, corresponding to parameters $\begin{bmatrix}
    1&1&1&1
\end{bmatrix}$, and do not depend on the starting guess of each optimization process.

The gradient-based approaches are stuck in the initial guess in most of the cases here analysed. In particular, in Tests 3a and 3b, where the starting point is the optimum obtained from the GA, both the gradient-based approaches do not improve the results obtained with the genetic. 
Only for Test 4b the gradient method converges to a point that does not coincide with the initial guess, but it does not improve the genetic results. Moreover, it reaches an efficiency value lower than the efficiency of the starting propeller of our pipeline.
We specify also that the values of the fitness reported in Table \ref{tab:grad} correspond to the \emph{penalized} efficiency of the propeller, i.e. the one expressed in \eqref{fitness_penalized}.

\RB{
\textbf{Remark:}\\
It is important to highlight that the genetic algorithm is characterized by a larger number of function evaluations with respect to gradient-based algorithms, resulting in a bigger computational time. The tuning of the hyperparameters, specified in Subsection \ref{subsec:genetic}, leads to a number of function evaluations equal to $180$ and $2750$, for the optimizations employing the \emph{standard} and \emph{fast} ROM, respectively.
Considering the fact that a function evaluation takes the same amount of time for gradient-based and genetic algorithms, the genetic computational time is \emph{one order of magnitude bigger} than the gradient-based optimization. In fact, Table \ref{tab:grad} shows that, in the constrained and \emph{fast} ROM optimization, the maximum number of function evaluations is $305$.
However, the gradient-based method produces results (in Table \ref{tab:grad}) that are often stuck in local minima and never approach the results of the genetic algorithm (reported in Table \ref{tab:genetic}).
Therefore, we can finally say that the genetic algorithm provides a robust and efficient method, that has been validated in many previous industrial shape optimization projects, such as ~\cite{DemoTezzeleMolaRozza2021JMSE, demo2021supervised, demoortaligustinrozzalavini2020bumi}.
}

\section{Conclusions}
\label{sec:conclusions}
The manuscript has illustrated a shape optimization pipeline exploiting data-driven ROM for improving the efficiency. The complexity of the here-discussed application, both in the geometry of the propeller and in the formulation of the original model, highlights the modularity and generality of the ROM based on data.

The parameterization of the blade\RB{, discussed in Section \ref{sec:param-def},} has been propagated to the mesh nodes, obtaining a series of discretized spaces that share the same topology, allowing us to apply the POD technique.
\RB{As presented in Subsection \ref{subsec:mesh-def},} the deformation applied to the mesh nodes maintains an overall quality similar to the original mesh and avoids reconstructing a new mesh for any new deformation. It must be said that due to the RBF technique applied in this phase, the computational cost of the mesh deformation is not negligible, but lower than the one required for constructing the mesh from scratch. This possible computational bottleneck will be however better addressed in future works in order to mitigate its impact over the entire pipeline.

\RB{A problem related to the computational burden of the mesh generation also arises in the online stage, when using a \emph{standard ROM} approach, as discussed in Subsection \ref{rom_1}.}
One possible solution --- that resolves the problem only at the online stage --- is explored here with the \emph{Fast ROM}\RB{, introduced in Subsection \ref{rom_2}}. The data-driven model is indeed built over the quadrature nodes of the blade faces, avoiding the deformation of the mesh since the objective function to minimize depends only on integral quantities over the blade surface.

\RB{As the reader can see in Subsection \ref{subsec:results-roms}, }the \rom{}s demonstrate in this work great precision, even thanks to the large offline simulation campaign performed to collect the initial high-fidelity snapshots. The computational cost required by the creation of this database is of course not negligible but has allowed a big reduction during the optimizations, allowing to test several methods and for the hyper-parameters tuning thanks to the huge gain in simulation time ($\SI[scientific-notation=false]{24}{\hour} = \SI{8.64e4}{\second}$ for the FOM against $\SI[scientific-notation=false]{15}{\second}$ for the Fast-ROM, more than $5000$ times faster). Future works will perform anyway the sensitivity analysis at varying the dimension of the solutions database, trying to provide guidelines for increasing accuracy in fixed computational budget contexts.

The efficiency of the reduced models allowed for testing different optimization procedures \RB{in Subsections \ref{sec:opt}}, comparing in this case the employment of methods based on the gradient of the objective functions \RB{in Subsection \ref{subsec:gradient}} and methods based on genetic strategies \RB{in Subsection \ref{subsec:genetic}}. Even if the gradient-based approach is the most rigorous one, its employment over a non-convex manifold (here computed by the RBF technique) produces poor results in terms of performance and shows its sensitivity to the selected starting point. The genetic optimization, requiring a larger number of evaluations, is able instead to produce better shapes in all the tests we performed.

\section{Tables}
\label{sec:tables}
\begin{table*}[htbp]
    \centering
    \caption{Accuracy analysis for the \textbf{standard} \rom{}. The accuracy is measured in terms of the mean of the \emph{k-fold cross validation} errors with $k=10$.}
    \begin{tabular}{p{0.3\textwidth}p{0.08\textwidth}p{0.08\textwidth}p{0.08\textwidth}}
    \toprule 
    \rom{} field &\rbf{} & \gpr{} & \knr{} \\
        \midrule
        Pressure & $0.776 \%$&  \cellcolor{green!10} $0.632\, \%$ & $9.651 \,\%$\\
        Wall shear stress (x)&  \cellcolor{green!10} $1.542 \,\%$ & $3.571 \,\%$ & $3.348 \,\%$ \\
        Wall shear stress (y)&  \cellcolor{green!10} $3.293 \,\%$ & $6.305 \,\%$ & $5.310 \,\%$ \\
        Wall shear stress (z)&  \cellcolor{green!10} $1.335 \,\%$ & $2.156 \,\%$ & $3.968 \,\%$ \\
        \bottomrule
    \end{tabular}
    \label{tab:rom_analysis_std}
\end{table*}

\begin{table*}[htbp]
    \centering
    \caption{Accuracy analysis for the \textbf{fast} \rom{}. The accuracy is measured in terms of the mean of the \emph{k-fold cross validation} errors with $k=10$.}
    \begin{tabular}{p{0.3\textwidth}p{0.08\textwidth}p{0.08\textwidth}p{0.08\textwidth}}
    \toprule 
     \rom{} field &\rbf{} & \gpr{} & \knr{} \\
        \midrule
        Pressure &  \cellcolor{green!10}$2.673 \%$& $5.857\, \%$ & $12.117 \,\%$\\
        Wall shear stress (x)& \cellcolor{green!10} $2.913 \,\%$ & $6.208 \,\%$ & $4.852 \,\%$ \\
        Wall shear stress (y)&  \cellcolor{green!10} $6.252\,\%$ & $12.233 \,\%$ & $7.489 \,\%$ \\
        Wall shear stress (z)& \cellcolor{green!10} $2.959 \,\%$ & $5.907 \,\%$ & $5.152 \,\%$ \\
        Normals (x)& \cellcolor{green!10} $4.688 \,\%$ & $11.312 \,\%$ & $4.984 \,\%$ \\
        Normals (y)&  $4.595 \,\%$ & $11.03 \,\%$ & \cellcolor{green!10} $4.487 \,\%$ \\
        Normals (z)& \cellcolor{green!10} $4.511 \,\%$ & $10.837 \,\%$ & $4.888 \,\%$ \\
        \bottomrule
    \end{tabular}
    \label{tab:rom_analysis_fast}
\end{table*}

\begin{table}[htbp]
    \centering
    \caption{Results of genetic optimization for the unconstrained and constrained cases, for both the standard and fast \rom{}s. $\Delta_{\text{ROM}}$ and $\Delta_{\text{FOM}}$ are the efficiency improvements with respect to the starting propeller of the optimization processes.}
    \begin{tabular}{p{0.1\textwidth}p{0.15\textwidth}p{0.07\textwidth}p{0.07\textwidth}p{0.07\textwidth}p{0.07\textwidth}p{0.08\textwidth}p{0.08\textwidth}}
    \toprule
    &  & \multicolumn{4}{c}{Parameters}&  &  \\
    \cmidrule{3-6}
    & Optimization & Def. \textbf{pitch} & Def. \textbf{camber} & Def. \textbf{chord length} & Def. \textbf{thickness} & $\Delta_{\rom{}}$  & $\Delta_{\fom{}}$ \\
    \midrule
         \multirow{2}{*}{\textbf{Standard}} & Unconstrained & $0.92$ & $1.03$ &  $0.70$ &  $0.75$ & $+5.15 \, \%$ & $+5.13 \, \%$\\
         &Constrained & $1.01$ & $0.81$ & $0.82$ & $0.99$ & $+0.79 \, \%$ & $+0.81 \, \%$ \\ 
         \midrule
         \multirow{2}{*}{\textbf{Fast}} & Unconstrained &  $0.91$ &  $0.89$ &  $0.74$ &  $0.70$ & $+3.76 \, \%$ & $+3.24 \, \%$ \\
         &Constrained & $0.99$ & $0.90$ & $0.77$ & $1.02$ & $+1.16 \, \%$ & $+0.80 \, \%$\\ 
         \bottomrule
    \end{tabular}
    \label{tab:genetic}
\end{table}

\begin{table}[htbp]
    \centering
    \caption{Results of gradient optimization for different initial guesses for the CG and L-BFGS-B methods. The results are reported in terms of number of function and gradient evaluations, initial and final points of the optimization, and fitness improvement with respect to the starting propeller.}
    \label{tab:grad}
    \begin{tabular}{M{0.02\linewidth}
    M{0.12\linewidth}
    M{0.22\linewidth}
    M{0.035\linewidth}
    M{0.035\linewidth}
    M{0.15\linewidth}
    M{0.22\linewidth}
    }
    \toprule
         Test & Method & Initial guess & Func. evals &  Grad. evals & $\Delta_{\text{ROM}}$ & Final point \\
         \midrule
         1a & CG & \multirow{2}{*}{$\begin{bmatrix} 1 & 1 & 1 & 1\end{bmatrix}$} &160 & 30 & \multirow{2}{*}{$+0 \%$} & \multirow{2}{*}{$\begin{bmatrix} 1 & 1 & 1 & 1\end{bmatrix}$}\\
         \cmidrule{1-2}
         \cmidrule{4-5}
         1b & L-BFGS-B &  &190 & 38 &  & \\
         \midrule
         2a & CG & \multirow{2}{*}{$\begin{bmatrix} 0.9 & 1.1 & 0.9 & 1.1\end{bmatrix}$ }&255 & 49 & \multirow{2}{*}{$-0.4 \%$} & \multirow{2}{*}{$\begin{bmatrix} 0.9 & 1.1 & 0.9 & 1.1\end{bmatrix}$}\\    
         \cmidrule{1-2}
         \cmidrule{4-5}
         2b & L-BFGS-B &  &305 & 61 & & \\
         \midrule
         3a & CG & \multirow{2}{*}{$\begin{bmatrix} 0.99 & 0.90 & 0.77 & 1.02\end{bmatrix}$} &247 & 48 & \multirow{2}{*}{$+1.1 \%$} & \multirow{2}{*}{$\begin{bmatrix} 0.99 & 0.90 & 0.77 & 1.02\end{bmatrix}$}\\      
         \cmidrule{1-2}
         \cmidrule{4-5}
         3b & L-BFGS-B &  &120 & 24 & & \\
         \midrule
         4a & CG & \multirow{2}{*}{$\begin{bmatrix} 0.9 & 0.8 & 0.8 & 0.9\end{bmatrix}$} & 151 & 28 & $-32.6 \%$& $\begin{bmatrix} 0.9 & 0.8 & 0.76 & 0.9\end{bmatrix}$\\  
         \cmidrule{1-2}
         \cmidrule{4-7}
         4b & L-BFGS-B & & 215 & 43 & $-0.1 \%$ & $\begin{bmatrix} 1 & 1 & 0.75 &  1.15\end{bmatrix}$\\
         \bottomrule
    \end{tabular}
    \end{table}

\section*{Acknowledgments}
We sincerely thank the industrial partners' Ing. Gianluca Gustin, Ing. Gianpiero Lavini, Ing. Edoardo Tagliamonte and Ing. Nicola Iona (from FINCANTIERI S.P.A.).
We also thank Francesco Andreuzzi and Gianmarco Gurioli for their preliminary work on this project.

This work was partially funded by INdAM-GNCS 2020-2021 projects, by European High-Performance Computing Joint Undertaking project Eflows4HPC GA N. 955558, by PRIN "Numerical Analysis for Full and Reduced
Order Methods for Partial Differential Equations" (NA-FROM-PDEs) project
by European Union Funding for
Research and Innovation --- Horizon 2020 Program --- in the framework
of European Research Council Executive Agency: H2020 ERC CoG 2015
AROMA-CFD project 681447 ``Advanced Reduced Order Methods with
Applications in Computational Fluid Dynamics'' P.I. Professor Gianluigi Rozza.

\subsection*{Author contributions}
Conceptualization: Anna Ivagnes, Nicola Demo; Methodology: Anna Ivagnes, Nicola Demo; Formal analysis and investigation: Anna Ivagnes; Writing - original draft preparation: Anna Ivagnes; Writing - review and editing: Nicola Demo; Funding acquisition: Gianluigi Rozza; Supervision: Gianluigi Rozza.

\newpage

\bibliographystyle{plain}
\bibliography{biblio}
\end{document}